\documentclass{elsarticle}
\usepackage{amssymb,amsmath}
\usepackage{graphicx,bbm,pstricks}
\usepackage{subfigure}
\usepackage{lineno}
\usepackage{bbm,multirow}
\newtheorem{definition}{Definition}[section]
\newcommand{\RARR}[3]{#1
  \;\displaystyle\mathop{\displaystyle\longrightarrow}^{#3}\; #2}

\newcommand{\LRARR}[4]{{\mbox{ \raise 0.4 mm \hbox{$#1$}}} \;
  \mathop{\stackrel{\displaystyle\longrightarrow}\longleftarrow}^{#3}_{#4}
  \; {\mbox{\raise 0.4 mm\hbox{$#2$}}}} 
\newcommand{\bX}{{\bf X}}

\newcommand{\bS}{{\bf S}}
\newcommand{\bF}{{\bf F}}
\newcommand{\bP}{{\bf P}}
\newcommand{\bnu}{{\boldsymbol{\nu}}}

\renewcommand{\red}[1]{#1}
\renewcommand{\blue}[1]{#1}
\begin{document}

\begin{frontmatter}
\nolinenumbers

\title{Constrained Approximation of Effective Generators for
  Multiscale Stochastic Reaction Networks and
  Application to Conditioned Path Sampling}
\journal{Journal of Computational Physics}
\author{Simon L. Cotter}
\address{School of Mathematics,
University of Manchester, Oxford Road,
Manchester, M13 9PL, United Kingdom;
e-mail: simon.cotter@manchester.ac.uk.
SC was funded by First Grant Award EP/L023989/1 from EPSRC.}

\begin{abstract}
Efficient analysis and simulation of multiscale \red{stochastic} systems of chemical
kinetics is an ongoing area for research, and is the source of many
theoretical and computational challenges. In this paper, we present a significant improvement to the constrained
approach, \red{which is a method for computing effective dynamics of
  slowly changing quantities in these systems, but which does not rely
  on the quasi-steady-state assumption (QSSA). The QSSA can cause
  errors in the estimation of effective dynamics for systems where the
  difference in timescales between the ``fast'' and ``slow'' variables
  is not so pronounced. 

This new application of the constrained approach} allows us to compute the effective generator of the
slow variables, without the need for expensive stochastic
simulations. This is \red{achieved by} finding the null space of the generator
of the constrained system. For complex systems where this is not
possible, or where the constrained subsystem is itself multiscale, the
constrained approach can then be applied \red{iteratively. This results
  in breaking the problem down into finding the solutions to many
  small eigenvalue problems, which can be efficiently solved using
  standard methods.} 

Since this methodology does not rely on the quasi
steady-state assumption, the effective dynamics that are
approximated are highly accurate, and in the case of systems with only
monomolecular reactions, are exact. We will demonstrate this with some
numerics, and also use the effective generators to sample paths of the
slow variables which
are conditioned on their endpoints\red{, a task which would be
computationally intractable for the generator of the full system.}
\end{abstract}

\end{frontmatter}
\section{Introduction}
Understanding of the biochemical reactions that govern cell function
and regulation is key to a whole range of biomedical and biological
applications and understanding mathematical modelling of gene regulatory networks has been an area
of huge expansion over the last half century. Due to the low copy
numbers of some chemical species within the cell, the random and
sporadic nature of individual reactions can play a key part in the
dynamics of the system, which cannot be well approximated by
ODEs\cite{Erban:2009:ASC}. Methods for the simulation of such a
system, such as
Gillespie's stochastic simulation algorithm
(SSA)\cite{gillespie1977exact}, \red{(or the similar
  Bortz-Kalos-Lebowitz algorithm\cite{bortz1975new} specifically for
  Ising spin systems),}
have 
been around for some decades. Versions which are more computationally efficient have also
been developed in the intermediate years\cite{Gibson:2000:EES,cao2004efficient}.

Unfortunately, their application to many systems can be very
computationally expensive, since the algorithms simulate every single
reaction individually. If the system is multiscale, i.e. there are some
reactions (fast reactions) which are happening many times on a timescale for which
others (slow reactions) are unlikely to happen at all, then in order for us to
understand the occurrences of the slow reactions, an unfeasible number
of fast reactions must be simulated. This is the motivation for
numerical methods which allow us to approximate the dynamics of the
slowly changing quantities in the system, without the need for
simulating all of the fast reactions.

For systems which are assumed to be well-mixed, there are many
different approaches and methods which have been developed. For
example the $\tau$-leap method\cite{gillespie2001approximate} speeds up the simulation by
timestepping by an increment within which several reactions may
occur. This can lead to problems when the copy numbers of one or more
of the species
approaches zero, and a number of different methods for overcoming this
have been presented\cite{tian2004binomial,auger2006r}.

Several other methods are based on the quasi steady-state assumption
(QSSA). This is the assumption that the fast variables converge in
distribution in a time which is negligible in comparison with the
rate of change of the slow variable. Through this assumption, a simple
analysis of the fast subsystem yields an approximation of the dynamics
of the slow variables. This fast subsystem can be analysed in several
ways, either through analysis and approximation\cite{cao2005slow}, or
through direct simulation of the fast subsystem\cite{weinan2005nested}.

Another approach is to approximate the system by a continuous
state-space stochastic differential equation (SDE), through the
chemical Langevin equation (CLE)\cite{gillespie2000chemical}. This
system can then be simulated using numerical methods for SDEs.
An alternative approach is to approximate only the slow variables by
an SDE. The SDE parameters
can be found using bursts of stochastic simulation of the system,
initialised at a particular point on the slow state space\cite{Erban:2006:GRN}, the so-called ``equation-free'' approach. This was
further developed into the constrained multiscale algorithm
(CMA)\cite{cotter2011constrained}, which
used a version of the SSA which also constrained the slow variables to
a particular value. Using a similar approach to \cite{cao2005slow}, the CMA
can similarly be adapted so that approximations of the invariant
distribution of this constrained system can be made without the need
for expensive stochastic simulations\cite{cucuringu2015adm}. However,
depending on the system, as with the slow-scale SSA, these approximations may
incur errors. \red{Work on how to efficiently approximate the results of multiscale
  kinetic Monte Carlo problems is also being undertaken in many different
  applications such as Ising models and lattice gas
  models\cite{novotny1995monte}.}

Analysis of mathematical models of gene regulatory networks (GRNs) is
important for a number of reasons. It can give us further insight into
how important biological processes within the cell, such as the
circadian clock\cite{vilar2002mechanisms}   or the cell cycle\cite{kar2009exploring} work. In order
for these models to be constructed, we need to observe how these
systems work in the first place. Many of the observation techniques, such as the DNA
microarray\cite{Schena20101995}, are notoriously subject to a large amount of
noise. Moreover, since the systems themselves are stochastic, the
problem of identifying the structure of the network from this data is
very difficult. As such, the inverse problem of
characterising a GRN from observations is a big challenge facing our
community\cite{golightly2014bayesian}.

One popular approach to dealing with inverse problems, is to use a
Bayesian framework. The Bayesian approach allows us to combine prior
knowledge about the system, complex models and the observations in a
mathematically rigorous way\cite{stuart2010inverse}. In the context of GRNs, we only have noisy observations of the
concentrations of species at a set of discrete times. As such, we have
a lot of missing information. This missing data can be added to the
state space of quantities that we wish to infer from the data that we
do have. This complex probability distribution on both the true
trajectories of the chemical concentrations, and on the network
itself, can be sampled from using Markov chain Monte Carlo (MCMC)
methods, in particular a Gibb's
sampler\cite{fearnhead2006exact}. Within this Gibb's sampler, we need
a method for sampling a continuous path for the chemical
concentrations given a guess at the reaction parameters, and our noisy
measurements. Exact methods for
sampling paths conditioned on their endpoints have been developed
\cite{fearnhead2006exact,rao2013fast}. 

\red{ In other applications, methods
for path analysis and path sampling have been developed, for example discrete path sampling databases for
  discrete time Markov chains\cite{trygubenko2006graph}, or where the
  probability of paths, rather than that of trajectories of discrete Markov processes can be
  used to 
  analyse  behaviour\cite{sun2006path}. In \cite{eidelson2012transition}, a method for
  transition path sampling is presented for protein folding, where the
  Markov chain has absorbing states. Other approaches for
  coarse-graining transition path sampling in protein folding also
  exist\cite{berezhkovskii2009reactive}. Other methods also exist for
  the simulation of rare events where we wish to sample paths
  transitioning from one stable region to another\cite{bolhuis20113}.}

The problems become even more difficult when, as is often the case,
the systems in question are also multiscale. This means that these
inverse problems require a degree of knowledge from a large number of
areas of mathematics. Even though many of the approaches that are
being developed are currently out of reach in terms of our current
computational capacity, this capacity is continually improving. In
this paper we aim to progress this methodology in a couple of areas.

\subsection{Conditioned path sampling methods}\label{sec:CPS}
We will briefly review the method presented in
\cite{fearnhead2006exact} for the exact sampling of conditioned paths
in stochastic chemical networks. Suppose that we have a Markov jump
process, possibly constructed from such a network, with a generator
$\mathcal{G}$. \red{The generator of such a process is the operator
$\mathcal{G}$ such that the master equation of the system can be
expressed as $$\frac{d{\bf p}}{dt} = \mathcal{G} {\bf p},$$ where
${\bf p}$ is the (often infinite dimensional) vector of probabilities
of being in a particular state in the system.} We wish to sample a path, conditioned on $X(t_0) = x_0$
and $X(t_1) = x_1$. Such a path can be found by creating a dominating
process (i.e. a process whose rate is greater than the fastest rate of
any transitions of the original system) with a uniform rate.

\begin{table}
\framebox{%
\hsize=0.97\hsize
\vbox{
\leftskip 6mm
\parindent -6mm
{\bf [1]} Define a dominating process to have transition rates given
by the matrix $\mathcal{M} = \frac{1}{\rho} \mathcal{G} + I$.

\smallskip

{\bf [2]} This process has uniformly distributed reaction events on the
time interval $[t_0,t_1]$. The number $r$ of such events is given by \eqref{eq:S:numR}.

\smallskip

{\bf [3]} Once $r = \hat{r}$ has been sampled, the type of each event
must be decided, by sampling from the distribution \eqref{eq:S:type},
starting with the first event. An event which corresponds to rate
$m_{i,i}$ indicates that no reaction event has occurred at this event.

\smallskip

{\bf [4]} Once all event types have been sampled, we have formed a
sample from the conditioned path space.

\par \vskip 0.8mm}
}
\caption{{\it A summary of the methodology presented in
    \cite{fearnhead2006exact}, for sampling paths of Markov-modulated
    Poisson processes, conditioned on their
  endpoints.}\label{CPS table}}
\end{table}

We define the rate to be greater than
the fastest rate of the process with generator $\mathcal{G}$, so
that $$\rho > \max_i \mathcal{G}_{i,i}.$$ 
Then we define the
transition operator of the dominant process by:
$$\mathcal{M} = \frac{1}{\rho} \mathcal{G} + I.$$
We can then derive the number of reaction events $N_U$ of the
dominating process in the time interval
$[t_0,t_1]$ by:
\begin{equation}\label{eq:S:numR}
\mathbb{P}(N_U = r) = \frac{\exp(-\rho t)(\rho t)^r/r! [\mathcal{M}^r]_{x_0,x_t}}{[\exp(\mathcal{G}t)]_{x_0,x_t}}.
\end{equation}
\red{Here the notation $[\cdot]_{a,b}$ denotes the entry in the matrix with
coordinates $(a,b) \in \mathbb{N}^2$.} A sample is taken from this distribution. The times $\{t^*_1, t^*_2,
\ldots t^*_{r}\}$ of all of the $r$
reaction events can then be sampled uniformly from the interval
$[t_0,t_1]$. The only thing that then remains is to ascertain which
reaction has occurred at each reaction event. This can be found by
computing, starting with $X(t_0) = x_0$, the probability distribution defined by:
\begin{equation}\label{eq:S:type}
\mathbb{P}(X(t^*_j)) = x | X(t^*_{j-1}) = x^*_{j-1}, X(t_1) = x_1) = \frac{
[\mathcal{M}]_{x^*_{j-1},x} [\mathcal{M}^{r-j}]_{x,x_1}}{[\mathcal{M}^{r-j+1}]_{x^*_{j-1},x_1}}.
\end{equation}
This method, summarised in Table \ref{CPS table}, exactly samples from the desired distribution, but
depending on the size and sparsity of the operator $\mathcal{G}$, it
can also be very expensive. In the context of multiscale systems with
a large number of possible states of the variables, the method quickly
becomes computationally intractable.

\subsection{Summary of Paper}
In Section \ref{Sec:CMA}, we introduce a version of the Constrained Multiscale
Algorithm (CMA), which allows us to approximate the effective
generator of the slow processes within a multiscale system. In
particular, we explore how stochastic simulations are not
required in order to compute a highly accurate effective generator. 
In
Section \ref{sec:comp}, we \red{consider the differences between the
  constrained approach, and the more commonly used quasi-steady state
  assumption (QSSA). In Section
\ref{sec:Nest}, we describe how the constrained approach can be
extended in an iterative nested structure for systems for whose constrained
subsystem is itself a large intractable multiscale system. By applying the
methodology in turn to the constrained systems arising from the
constrained approach, we can make the
analysis of
highly complex and high dimensional systems computationally
tractable. In Section
\ref{sec:num}, we present some analytical and numerical results, aimed
at presenting the advantages of the CMA over other approaches. This includes some
examples of conditioned path sampling using effective generators
approximated using the CMA. Finally, we will summarise our findings in
Section \ref{sec:conc}.}
\section{The Constrained Multiscale Algorithm}\label{Sec:CMA}
The Constrained Multiscale Algorithm was originally designed as a multiscale
method which allowed us to compute the effective drift and diffusion
parameters of a diffusion approximation of the slow variables in a
multiscale stochastic chemical network. The idea was simply to
constrain the original dynamics to a particular value of the slow
variable. This can be done through a simple alteration of the original
SSA by Gillespie\cite{gillespie1977exact}. \red{First, a (not
  necessarily orthogonal) basis is
  found for the system in terms of ``slow'' and ``fast''
  variables, $[{\bf S} = [S_1,S_2,\ldots], {\bf F} =
  [F_1,F_2,\ldots]]$. Slow variables are not affected by the most
  frequently firing reactions in the system. Then, }as shown in
\cite{cotter2011constrained}, the SSA is computed as normal, until one
of the slow reactions \blue{(a reaction which alters the value of the slow variable(s))} occurs. After the reaction has occurred, the
slow variable is then reset to its original value, in such a way that
the fast variables are not affected. \blue{This is equivalent to projecting
the state of the system, after each reaction, back to the desired
value of the slow variable, whilst also
preserving the value(s) of the fast variable(s).} The constrained SSA is given in
Table \ref{CSSA table}. \blue{Here the $\alpha_i({\bf X}(t))$
  denote the propensity of the reaction $R_i$ when the system is in
  state $\bX(t) = [X_1(t),X_2(t),\ldots]$, where $\Delta t
  \alpha_i({\bf X}(t)$ is the probability that
  this reaction will fire in the infinitesimally small time interval
  $(t, t+\Delta t)$ with $1 \gg \Delta t > 0$. The stoichiometric
  vectors $\bnu_i$ denote the change in the state vector $\bX(t)$ due
  to reaction $R_i$ firing.}

\blue{In order to describe the constrained approach, we first
  introduce some definitions that will be helpful.}

\blue{\begin{definition}{Constrained Projector:} Given a basis of the state space $\bX =
    [X_1,X_2,\ldots,X_N]$ with $N_f$
    fast variables ${\bf F} = [F_1, F_2, \ldots, F_{N_f}]$ and $N_s$
    slow variables ${\bf S} = [S_1, S_2, \ldots, S_{N_s}]$, the
    constrained projector $\mathcal{P}_{\bf S}:\mathbb{N}_0^N \to
    \mathbb{N}_0^N$ for a given value of ${\bf S}$ preserves the
    values of the fast variables, whilst mapping the values of the
    slow variables to ${\bf S}$:
\begin{equation}
\mathcal{P}_{\bf S} ([\hat{\bf S},\hat{\bf F}]) = [\bS,\hat{\bF}]
\qquad \qquad \forall  ([\hat{\bf S},\hat{\bf F}]) \in \mathbb{N}_0^N.
\end{equation}
\end{definition}}

\blue{\begin{definition}{Constrained Stoichiometric Projector:}
Given a basis of the state space $\bX =
    [X_1,X_2,\ldots,X_N]$ with $N_f$
    fast variables ${\bf F} = [F_1, F_2, \ldots, F_{N_f}]$ and $N_s$
    slow variables ${\bf S} = [S_1, S_2, \ldots, S_{N_s}]$, the
    constrained stoichiometric projector $\mathcal{P}:\mathbb{N}_0^N \to
    \mathbb{N}_0^N$ maps any non-zero elements of the slow coordinates
    to zero, whilst preserving the values of the fast coordinates:
\begin{equation}
\mathcal{P}([{\bf S},{\bf F}]) = [{\bf 0},{\bF}]
\qquad \qquad \forall  ([{\bf S},{\bf F}]) \in \mathbb{N}_0^N.
\end{equation}\label{def:CSP}
\end{definition}}

\blue{\begin{definition}{Constrained Subsystem:}
Given a system with $N_R$ reactions $R_1, R_2, \ldots,
R_{N_R}$ with propensity functions $\alpha_i(\bX)$ and stoichiometric
vectors ${\boldsymbol{\nu}}_i \in \mathbb{N}_0^N$, the constrained subsystem is the system that arises from applying the
constrained projector $\mathcal{P}_{\bS}$ to the state vector after
each reaction in the system. This is equivalent to applying the
constrained stoichiometric projector $\mathcal{P}$ to each of the stoichiometric
vectors in the system. This may leave some reactions with a null
stoichiometric vector, and so these reactions can be removed from the
system. This projection can lead to aphysical systems where one or more variables
may become negative; in these cases we set the propensities of the
offending reactions at states where a move to a negative rate is
possible, to zero.
\end{definition}}

\blue{
\begin{table}
\framebox{
\hsize=0.97\hsize
\vbox{
\leftskip 6mm
\parindent -6mm 

{\bf [1]} Define a basis of the state space in terms of slow and fast
variables.

\smallskip

{\bf [2]} Initialise the value of the state, $\bX(t_0) = {\bf x}$.

\smallskip

{\bf [3]} Calculate propensity functions at the current state $\alpha_i({\bX}(t))$.

\smallskip

{\bf [4]} Sample the waiting time to the next reaction in the
system
\begin{equation*}
  \tau = -\frac{\log\left(u\right )}{\alpha_0({\mathbf X}(t))},
  \quad \mbox{where} \quad \alpha_0({\mathbf X}(t)) 
  = \sum_{k=1}^M \alpha_k({\mathbf X}(t)), \blue{\qquad u \sim U([0,1]).}
\end{equation*}

\smallskip

{\bf [5]} Choose one $j\in\{1,\ldots,M\}$, with probability 
$\alpha_j/\alpha_0$, and perform reaction $R_j$, with stoichiometry
which has been projected using the constrained stoichiometric
projector:
\begin{equation*}
\bX(t+\tau) = \bX(t) + \mathcal{P}(\nu_j).
\end{equation*}

\smallskip

{\bf [6]} Repeat from step {\bf [3]}.
}
}
\caption{{\it The \emph{Constrained} Stochastic Simulation Algorithm
    (CSSA) using the constrained stoichiometric projector given in
    Definition \ref{def:CSP}. Simulation starts with $S=s$ where $s$ is a given value of
  the slow variable.}\label{CSSA table}}
\end{table}

}

Let us illustrate this using an example which we shall be using
later in the paper.
\begin{eqnarray}
R_1 &:& \qquad\qquad \emptyset \, \overset{k_1}{\longrightarrow} \, X_1,
\nonumber \\
R_2 &:& \qquad\qquad X_2 \, \overset{k_2}{\longrightarrow} \, \emptyset,
\label{eq:lin} \\
R_3 &:& \qquad\qquad X_1 \, \overset{k_3}{\longrightarrow} \, X_2,
\nonumber \\
R_4 &:& \qquad\qquad X_2 \, \overset{k_4}{\longrightarrow} \, X_1.  
\nonumber 
\end{eqnarray}
In certain parameter regimes,
this system is multiscale, with reactions $R_3$ and $R_4$ occurring
many times on a time scale for which reactions $R_1$ and $R_2$ are
unlikely to happen at all. The variable $S = X_1 + X_2$ is unaffected
by these fast reactions, and as such is a good candidate for the slow
variable which we wish to analyse. \red{A discussion about how the
  fast and slow variables could be identified is given in Section \ref{sec:conc}.} We have two choices for the fast
variable, either $F = X_1$ or $F=X_2$, \blue{in order to form a basis of the
state space along with the slow variable $S$.} As detailed in
\cite{cotter2011constrained}, it is preferable (although not
essential) to pick fast variables that are not involved in zeroth order
reactions. Therefore, in this case, we choose $F=X_2$. \blue{Following
  the projection of the stoichiometric vectors using the constrained
  projector, the
constrained system can be written in the following way:}
\begin{eqnarray}
C_1 &:& \qquad \qquad X_1 + X_2 = S,\nonumber\\
R_2 &:& \qquad\qquad X_2 \, \overset{k_2}{\longrightarrow} \, X_1, \nonumber
\\
R_3 &:& \qquad\qquad X_1 \, \overset{k_3}{\longrightarrow} \, X_2, 
\label{eq:lin:con} \\
R_4 &:& \qquad\qquad X_2 \, \overset{k_4}{\longrightarrow} \, X_1.  
\nonumber 
\end{eqnarray}
Note that reaction $R_1$ has disappeared completely, since \blue{only
  involves changes to the slow variable, and as such after projection,
  the stoichiometric vector is null, and the reaction can be removed.
The stoichiometry of reaction $R_2$ has been altered as it involves a
change in the slow variables. If this
reaction occurs, the slow variable is reduced by one. We are not
permitted to change the fast variable $X_2$ in order to
reset the slow variable to its original value, and therefore we must
increase $X_1$ by one, giving us a new stoichiometry for this reaction.}

In the original CMA, statistics were taken regarding the frequency of
the slow reactions, at each point of the slow domain, and were used to
construct the effective drift and diffusion parameters of an effective
diffusion\cite{cotter2011constrained,errorpaper} process. However, this constrained approach can also be used to
compute an effective generator for the discrete slow process, as we will now demonstrate.
The CMA can be very costly, due to the large computational burden of
the stochastic simulations of the constrained system. In this section,
we will introduce a method for avoiding the need for these
simulations, whilst also significantly improving accuracy.

The constrained systems can often have a very small state space (which
we will denote $\Gamma(s)$), since
they are constrained to a single value of the slow variables. For
example, for the constrained system \eqref{eq:lin:con}, there are only
$\left \lfloor \frac{S}{2} \right \rfloor$ possible states. Such a
system can easily be fully analysed. For example, the invariant
distribution can be found by characterising the one-dimensional
null space of the generator matrix of the constrained process. For small to medium-sized systems,
this is far more efficient than exhaustive Monte Carlo
simulations. For other systems with larger constrained state spaces,
stochastic simulation may still be the best option, although in
Section \ref{sec:Nest} we show how the constrained approach can be
applied iteratively until the constrained subsystem is easily analysed.

Suppose that we have a constrained system with $N_f$ fast variables,
$F_1, F_2, \ldots, F_{N_f}$. The generator for the constrained system
with $S=s$
is given by $\mathcal{G}_F(s)$. Since the system is ergodic, there is a
one-dimensional null space for this generator. This can be found by
using standard methods for identifying eigenvectors, by searching for
the eigenvector corresponding to the eigenvalue equal to zero. Krylov
subspace methods allow us to find these eigenvectors with very few iterations. Suppose
we have found such a vector \blue{${\bf v} = [v_1,v_2,\ldots]$,} such that
$$\mathcal{G}_F(s) {\bf v} = 0.$$
Then our approximation to the invariant distribution of this system is given by the
discrete probability distribution represented by the vector
\blue{$${\bf p}(s) = [p_1(s),p_2(s),\ldots] = \frac{\bf{v}}{\sum v_i}.$$}
Our aim is now to use this distribution to find the effective
propensities of the slow reactions of the original system.

Suppose that we have $N_s$ slow reactions in the original system. Each
has an associated propensity function $\alpha_1(S,F), \alpha_2(S,F),
\ldots, \alpha_{N_s}(S,F).$ We now simply want to find the expectation
of each of these propensity functions with respect to the probability
distribution ${\bf p}(s)$:\blue{
\begin{eqnarray}\label{effprop}
\mathbb{E} ( \alpha_i(S, \cdot)) = \sum_i p_i(s)\alpha_i(S,f).
\end{eqnarray}}
Having computed this expectation for all of the slow propensities,
over all required values of the slow variable, then an effective
generator for the slow variable can be constructed.

\begin{table}
\framebox{%
\hsize=0.97\hsize
\vbox{
\leftskip 6mm
\parindent -6mm
{\bf [1]} For each value of the slow variable $S=s \in \Omega$, compute the
generator $\mathcal{G}_s$ of the constrained subsystem.

\smallskip

{\bf [2]} Find the zero eigenvector \blue{${\bf v} = [v_1,v_2,\ldots]$} of $\mathcal{G}_s$, and
let ${\bf p}(s) = \frac{\bf{v}}{\sum v_i}$.

\smallskip

{\bf [3]} Approximate the effective propensities at each point $s \in
\Omega$ using \eqref{effprop}.

\smallskip

{\bf [4]} Construct an effective generator $\mathcal{G}$ of the slow
processes of the system using these effective propensities.

\par \vskip 0.8mm}
}
\caption{{\it The CMA approach to approximating the effective
    generator $\mathcal{G}$ of the slow variables \blue{on the
      (possibly truncated) domain $S \in \Omega$,} without the need
    for stochastic simulations.}\label{CMA table}}
\end{table}

\section{Comparing the CMA and QSSA approaches}\label{sec:comp}
A very common approach to approximating the dynamics of slowly
changing quantities in multiscale systems, is to invoke the quasi
steady-state assumption (QSSA). The assumption is that the fast and
slow variables are operating on sufficiently different time scales
that it can be assumed that the fast subsystem enters equilibrium
instantaneously following a change in the slow variables, \red{and
therefore is unaffected by the slow reactions.} This
assumption means that if the fast subsystem's invariant distribution
can be found (or approximated), then the effective
propensities of the slow reactions can be computed. However,
as demonstrated in \cite{errorpaper}, this assumption incurs an
error, and for systems which do not have a large difference in time
scales between the fast and slow variables, this error can be
significant.

The CMA does not rely on the QSSA, and is able to take
into account the effect that the slow reactions have on the invariant
distribution of the fast variables, conditioned on a value of the slow
variables. In a true fast-slow system, this will yield the same
results as the QSSA, but for most systems of interest, the constrained
approach will
have a significant increase in accuracy. If we follow the
approach outlined in Table \ref{CMA table}, we don't even need to
conduct any stochastic simulations to approximate the effective
dynamics.

The assumptions for the CMA are weaker than the
QSSA, namely that we assume that the dynamics of the slow variable(s) can
be approximated by a Markov-modulated Poisson process, independently
of the value of the fast variables. This means that we have made
the assumption that the current value of the fast variables has no effect
on the transition rates of the slow variables once a slow reaction has
occurred. This is subtly weaker than the QSSA, and importantly the
effect of the slow reactions on the invariant distribution of the fast
variables is accounted for. \red{Note that this may necessitate a slow
variable which has more than one dimension, for example in oscillating
systems for which the effective dynamics cannot be approximated by a
one dimensional Markov process. Consideration of such systems is an
area for future work.}

\section{The Nested CMA}\label{sec:Nest}
There will be many systems for which the constrained subsystem is
itself a highly complex and multiscale system. In this event, it will
not be feasible to find the null space of a sensibly truncated
generator for the constrained subsystem. Therefore, we need to
consider how we might go about approximating this. Fortunately, we
already have the tools to do this, since we can iteratively apply the
CMA methodology to this subsystem. This is analogous to the nested
strategy proposed in the QSSA-based nested
SSA\cite{weinan2005nested}. 

This nested approach allows us to reduce much more complex systems in
an accurate, computationally tractable way. The problem of finding the
null space of the first constrained subsystem is divided into finding
the null space of many small generators, through further constraining. An example of this nested approach will be demonstrated in Section \ref{sec:num:nest}.

\red{\section{Examples}\label{sec:num}
In this section we will present some analytical and numerical results produced using
the CMA approach for three different examples.} \blue{In order to give
an indication of the computational cost of the algorithms, we include
the runtime of certain operations. All numerics were performed using
MATLAB on a mid-2014 MacBook Pro. Disclaimer: the implementations used
are not highly
optimised, and these runtimes are purely given as an indication of the
true costs of a well implemented version.
}
\subsection{A Simple Linear System}\label{sec:lin}
\red{First we
consider a simple linear system, in order to demonstrate that the CMA
approximation of the effective generator of the slow variable is exact
in the case of systems with only monomolecular reactions, which is
 in contrast to the approximation found using a more standard QSSA-based approach. Let us illustrate this by returning to the example given by the linear
system
\eqref{eq:lin}, first analysing it using the QSSA.
\begin{eqnarray*}
R_1 &:& \qquad\qquad \emptyset \, \overset{k_1}{\longrightarrow} \, X_1,
\nonumber \\
R_2 &:& \qquad\qquad X_2 \, \overset{k_2}{\longrightarrow} \, \emptyset,
\\
R_3 &:& \qquad\qquad X_1 \, \overset{k_3}{\longrightarrow} \, X_2,
\nonumber \\
R_4 &:& \qquad\qquad X_2 \, \overset{k_4}{\longrightarrow} \, X_1.  
\nonumber 
\end{eqnarray*}
We will consider this system in the following parameter setting:
\begin{equation} \label{eq:lin:params}
k_1V = 20, \qquad k_2 = 1, \qquad k_3 = 5, \qquad k_4 = 5.
\end{equation}}
\blue{Here $V$ denotes the volume of the well-mixed thermally-equilibrated reactor.}
\subsubsection{QSSA-based analysis} \label{sec:lin:qssa}
The QSSA tells us that the fast subsystem
(made up of reactions $R_3$ and $R_4$) reaches probabilistic
equilibrium on a timescale which is negligible in comparison with the
timescale on which the slow reactions are occurring. Therefore we may
treat this subsystem in isolation with fixed $S$:

\begin{equation}
\LRARR{X_1}{X_2}{k_3}{k_4}, \qquad S = X_1 + X_2.\nonumber
\end{equation}

This is a very simple autocatalytic reaction system, for which a great
deal of analytical results are available. For instance, we can compute
the invariant distribution for this system\cite{jahnke2007solving},
which gives us that $X_2$ is a binomial random variable
$$X_2 \sim \mathcal{B} \left (\cdot, S, \frac{k_3}{k_3 + k_4} \right
).$$
Therefore, we can compute the conditional expectation
$\mathbb{E}(X_2|S) = \frac{k_3S}{k_3 + k_4} $ in this 
fast subsystem, and use this to approximate the effective rate of
reaction $R_2$. Therefore, the effective slow system is given by the
reactions:
\begin{equation}\label{eq:lin:eff}
\emptyset \, \overset{\hat{k}_1}{\longrightarrow} \, S  \, \overset{\hat{k}_2}{\longrightarrow} \, \emptyset,
\end{equation}
where
\begin{equation}
\hat{k}_1 = k_1, \qquad \hat{k}_2 = \frac{k_2\mathbb{E}(X_2)}{S} = \frac{k_2k_3}{k_3 + k_4}.\nonumber
\end{equation}
We can compute the invariant distribution for this effective
system\cite{jahnke2007solving}, which in this instance is a Poisson distribution:
\begin{equation}\label{eq:lin:QSS:ID}
S \sim \mathcal{P}\left (\frac{k_1V(k_3 + k_4)}{k_2k_3} \right ).
\end{equation}
We can quantify the error we have made in using the quasi-steady state
assumption by, for example, comparing this distribution with the true
invariant distribution. Once again, using the results of
\cite{jahnke2007solving}, we can compute the invariant distribution of
the system \eqref{eq:lin}, which is a multivariate Poisson
distribution:
\begin{equation} 
[X_1,X_2] \sim \mathcal{P}(\bar{\lambda}_1,\bar{\lambda}_2), \nonumber
\end{equation}
where $\bar{\lambda}_1 = \frac{k_1V(k_2 + k_4)}{k_2k_3}$, and
$\bar{\lambda}_2 = \frac{k_1V}{k_2}$. Trivially one can compute the
marginal distribution on the slow variable $S$:
\begin{eqnarray*}
\mathbb{P}(S=s) &=& \sum_{n=0}^s
\frac{\bar{\lambda}_1^n}{n!}\frac{\bar{\lambda}_2^{s-n}}{(s-n)!}\exp(-(\bar{\lambda}_1
+ \bar{\lambda}_2)),\\
&=& \frac{(\bar{\lambda}_1
 + \bar{\lambda}_2)^s}{s!}\exp(-(\bar{\lambda}_1
+ \bar{\lambda}_2)).
\end{eqnarray*}
Therefore $S$ is also a Poisson variable with intensity $\lambda = \bar{\lambda}_1
+ \bar{\lambda}_2 = \frac{k_1V(k_2+k_3 + k_4)}{k_2k_3}$, which differs
from the intensity approximated invariant density
\eqref{eq:lin:QSS:ID} by $\frac{k_1V}{k_3}$. Note that $k_3$ is one of
the fast rates, and $k_1V$ is one of the slow rates, and therefore as
the difference in timescales of the fast and slow reactions increases,
this error decreases to zero, so that the QSSA gives us an
asymptotically exact approximation of the slow dynamics. \red{For systems
with a finite timescale gap, the QSSA approximation
will incur error over and above the error incurred in any approximation of the marginalised
slow process by a Markov process.}

\subsubsection{CMA analysis}\label{sec:lin:CMA}
For comparison, let us compute approximations of the effective slow
rates by using the CMA. The CMA for this system tells us that we need
to analyse the constrained system \eqref{eq:lin:con}.
\begin{eqnarray*}
C_1 &:& \qquad \qquad X_1 + X_2 = S,\nonumber\\
R_2 &:& \qquad\qquad X_2 \, \overset{k_2}{\longrightarrow} \, X_1, \nonumber
\\
R_3 &:& \qquad\qquad X_1 \, \overset{k_3}{\longrightarrow} \, X_2, 
\\
R_4 &:& \qquad\qquad X_2 \, \overset{k_4}{\longrightarrow} \, X_1.  
\nonumber 
\end{eqnarray*}
The constrained system in this example only
contains monomolecular reactions, and as such can be analysed using
the results of \cite{jahnke2007solving}. The invariant distribution
for this system is a binomial, such that
$$X_2 \sim \mathcal{B} \left (\cdot, S, \frac{k_3}{k_2+k_3 + k_4} \right
).$$
Using this, we can compute the effective propensity of reaction $R_2$,
$$\bar{\alpha}_2(S) = k_2\mathbb{E}(X_2|S) =\frac{k_2k_3S}{k_2+k_3 +
  k_4},$$
giving us the effective rate $\bar{k_2} = \frac{k_2k_3}{k_2+k_3 +
  k_4}$. The invariant distribution of \eqref{eq:lin:eff} with this
effective rate for $\bar{k_2}$ is once again a Poisson distribution
with intensity $$\lambda= \frac{k_1V(k_2+k_3 + k_4)}{k_2k_3},$$
which is \emph{identical} to the intensity of the true distribution on
the slow variables. In other words, for this example, the CMA produces
an  approximation of the effective dynamics of the slow
variables for this system, whose invariant distribution is identical
to the marginal invariant distribution of the slow variables in the
full system.  The constrained approach corrects for the effect of the slow reactions
on the invariant distribution of the fast variables. In this and other
examples of systems with monomolecular reactions, the constrained
approach gives us a system whose invariant distribution is exactly
equal to the marginal distribution on the slow variables for the full
system. Another example is presented in Section \ref{sec:num:nest},
for which the constrained system is itself \red{too large to easily compute
expectations directly through its generator}, and requires
another iteration of the CMA to be applied.

For this example, we did not even
need to compute the invariant distributions of the constrained systems
numerically. In Section \ref{sec:bis}, we will come
across a system for which it is necessary to numerically compute the
invariant distribution of the constrained system.

\red{
\subsubsection{Comparison of approximation of invariant densities}

\begin{figure}[htp] 
\begin{center} 
\includegraphics[width=0.6\textwidth]{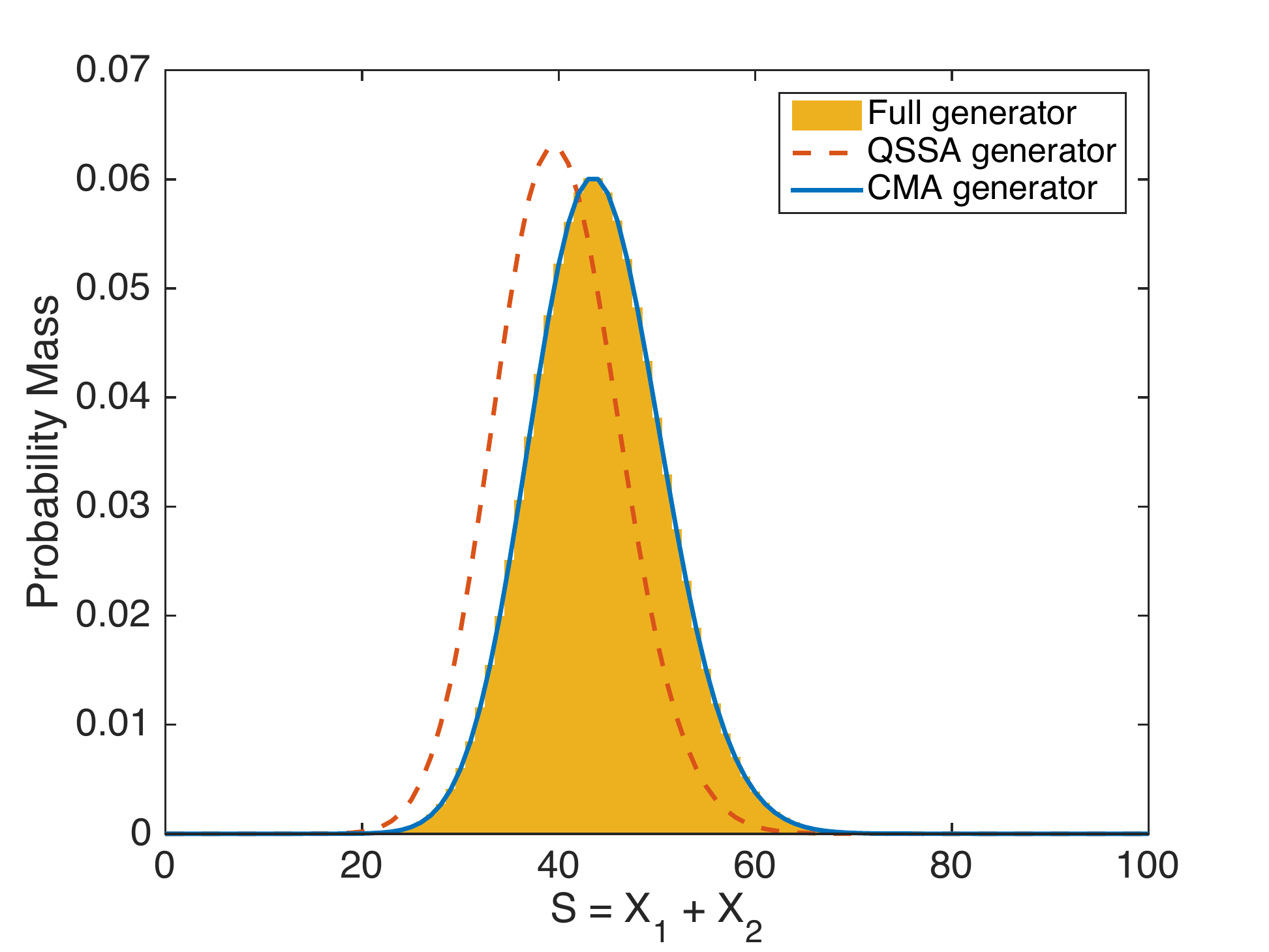}
\caption{\it Approximations of the invariant
  distribution of the slow variable $S=X_1 + X_2$ of system
  \eqref{eq:lin} with parameters \eqref{eq:lin:params} through
  marginalisation of the distribution of the full system (histogram), of the effective generator computed using the
  CMA (solid line) and computed using the QSSA (dashed line).\label{Fig:simple_pdfs}}
\end{center} 
\end{figure} 

Figure \ref{Fig:simple_pdfs} shows the invariant distributions of the
slow variables $S=X_1 + X_2$ in the parameter regime
\eqref{eq:lin:params}, computed by marginalising the invariant
distribution of the full system, and from the CMA and QSSA as outlined
above. The CMA exactly matches the true distribution, as both are
Poisson distributions with rate $\lambda = 44$. The QSSA incorrectly
approximates the effective rate of $R_4$, and as such is a Poisson
distribution with rate $\lambda = 40$. The relative error of the CMA
for this problem is zero, and for the QSSA is $4.322 \times
10^{-1}$. 


\subsubsection{Conditioned path sampling using effective generators}}
The approaches described in Section \ref{sec:CPS} hit problems when the system for which we are trying
to generate a conditioned path is multiscale. In a multiscale system,
the rate $\rho$ of the dominating process will be very large, and as
such the number of reaction events will be large, even if the path we
are trying to sample is short. Therefore $M^r$ is likely to be a full
matrix, and the number of calculations of \eqref{eq:S:type} will be
large. Moreover, the size of $M$ is also likely to be large, since for
each value $S=s$ of the slow variable, there are many states, one for
each possible value of the fast variable. All of
these factors make the problem of computing a conditioned path in such
a scenario computationally intractable.

Considering once more the system \eqref{eq:lin}, naturally we
cannot store the actual generator of this system, since the system is
open and as such the generator is an infinite dimensional
operator. However, the state space can be truncated carefully in such
a way that the vast majority of the states with non-negligible
invariant density are included, but an infinite number of highly
unlikely states are presumed to have probability zero. Note
  that this means that we are effectively sampling paths
  satisfying $S(t_0) = s_1$, $S(t_1) = s_2$ conditioned on $S(t) \in
  \Omega \quad \forall t$. However, even
with careful truncation the number of states can be prohibitively
large.

Suppose that we truncate the domain for this system
to $$\Omega = \left \{ [X_1,X_2] | X_1,X_2 \in \{0,1,\ldots,200 \right
\}.$$
This truncated system has $201^2 = 40401$ different states, and
therefore the generator $\mathcal{G} \in \mathbb{R}^{40401 \times
  40401}$. Although this matrix is sparse, the matrix exponential
required in \eqref{eq:S:numR} is full, as is $M^r$ for moderate $r
\in \mathbb{N}$. A full matrix of this size stored at double precision
would require over 13GB of memory. So even for this system, the most
simple multiscale system that one could consider, the problem of
sampling conditioned paths is computationally intractable.

In comparison, suppose that we use a multiscale method such as the CMA
to approximate the effective rates of the slow reactions. Then, for
the same $\Omega$, we only have 401 possible states of the slow
variable, a reduction of $99.25\%$. The effective generator
$\mathcal{G} \in \mathbb{R}^{401\times 401}$ would then only require
1.29MB to be stored as a full matrix in double precision. The
dominating process for this system must now have rate \red{$\rho>201.4$},
instead of \red{$\rho>1220$}, which is over \red{6} times bigger. This means
far fewer calculations of \eqref{eq:S:type}. What is
more, as we saw in Section \ref{sec:lin:CMA}, for some systems the CMA
\emph{exactly} computes the effective dynamics of the slow variables,
with no errors. 

\red{The system \eqref{eq:lin}, in order to highlight more effectively
the differences between the CMA and a QSSA-based approach, is in a
parameter range \eqref{eq:lin:params}, for which the
difference in time scales between the ``fast'' and ``slow'' variables
is relatively small, and of course for systems with larger timescale
difference, the difference in $\rho$ between the full and effective
generators would be far larger.}

Naturally, this approach only allows us to sample the paths of the
slow variables. However, the fast variables, if required, can easily
be sampled after the fact, using an adapted Gillespie approach which
samples the fast variables given a trajectory of the slow variables.

As we have just demonstrated in the previous section, the CMA can be used to compute an effective generator
for the slow variable $S = X_1 + X_2$ in the system \eqref{eq:lin}, with
parameters \eqref{eq:lin:params}, whose invariant distribution is
exactly that of the slow variable in the full system without the
multiscale reduction. Moreover, this can be achieved with no Monte
Carlo simulations, since the constrained subsystem contains only
monomolecular reactions, and as such its invariant distribution can be
exactly computed\cite{jahnke2007solving}. 

At this juncture, we simply need to apply the method of Fearnhead and
Sherlock\cite{fearnhead2006exact} to the computed effective generator in order to be able sample paths
conditioned on their endpoints. Suppose we wish to sample paths
conditioned on $S(t_0 = 0) = \red{44} = S(t_1 = 10)$. The invariant
distribution of this system, as shown previously in this paper, is a
Poisson distribution with mean \red{$\lambda= \frac{k_1V(k_2+k_3 +
  k_4)}{k_2k_3} = 44$. Therefore, we are attempting to sample paths
which start and finish at the the mean of the invariant distribution,
which in itself is not a particularly interesting thing to do, but it
will allow us to highlight again the advantages of using the CMA over
QSSA-based approaches.}

Since the system is open, we are required to truncate the domain in
order to be able to store and manipulate the effective generator. We
truncate the domain to $\Omega = \{[X_1,X_2] | S = X_1 + X_2 \leq
400\}$. Therefore we aim to sample paths $$\large \{S(t),\, t\in [0,10] \,|\,  S(0) =
 \red{44} = S(10),\, S(t) \in \Omega \, \forall t\in [0,10] \large\}.$$

As the number of possible states of the slow variable is relatively
small, it was possible to compute and store full matrices for $M^r$ as
required in \eqref{eq:S:numR} and \eqref{eq:S:type} for \red{$r \in
1,2,\ldots,2369$. $r$ has an upper bound of 2369 as the cumulative
mass function for the probability distribution \eqref{eq:S:numR} is
within machine precision of one at $r=2369$.} Storing all powers of the
matrices is clearly not the most efficient way to
implement this algorithm, but for this example was possible without
any intensive computations, and with minimal numerical error. We will
present a more efficient approach in the next section.

\red{

\begin{figure}[htp] 
\begin{center} 
\mbox{ 
\subfigure[]{\scalebox{0.3}{\includegraphics{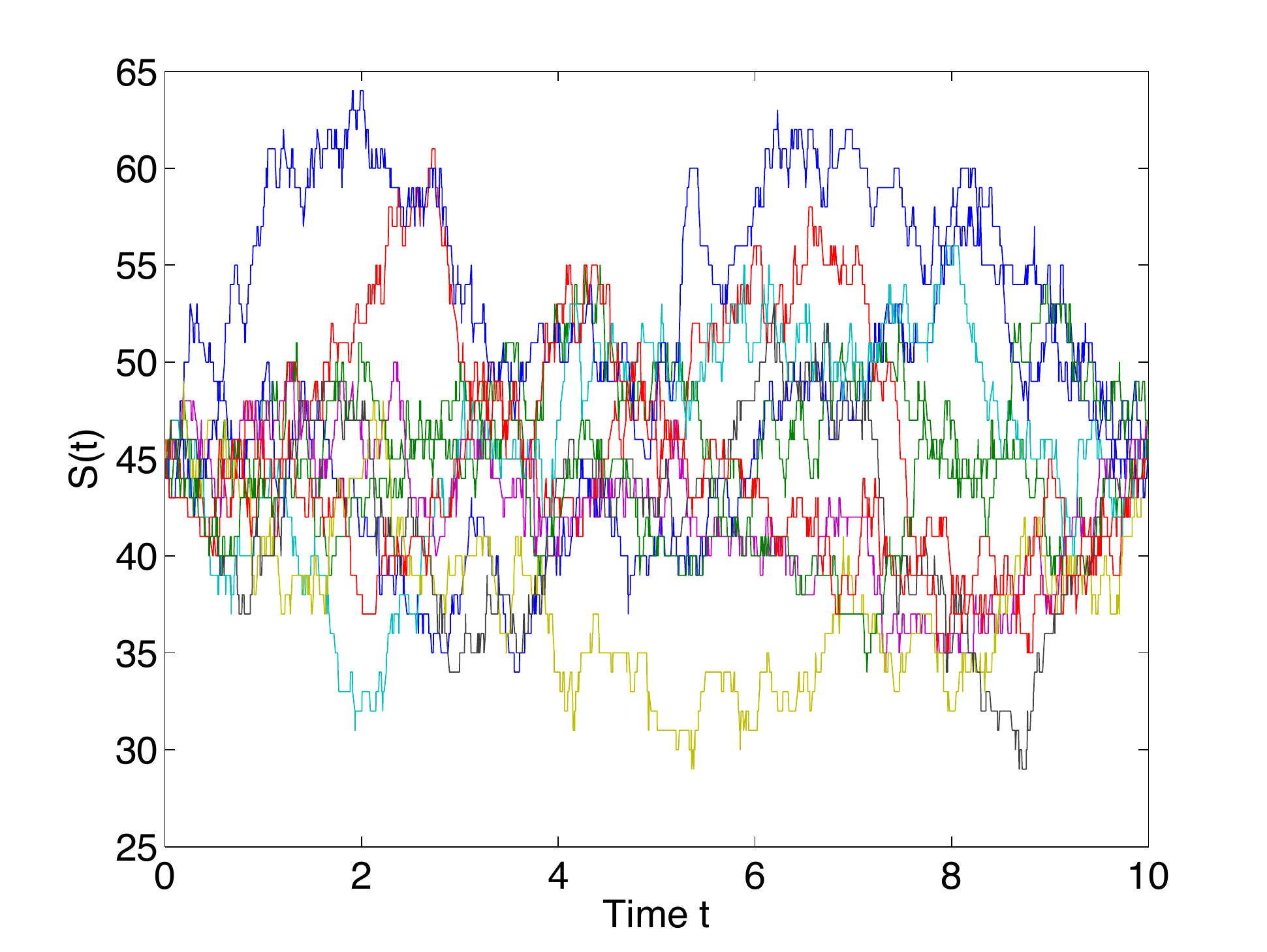} }}
\subfigure[]{\scalebox{0.3}{\includegraphics{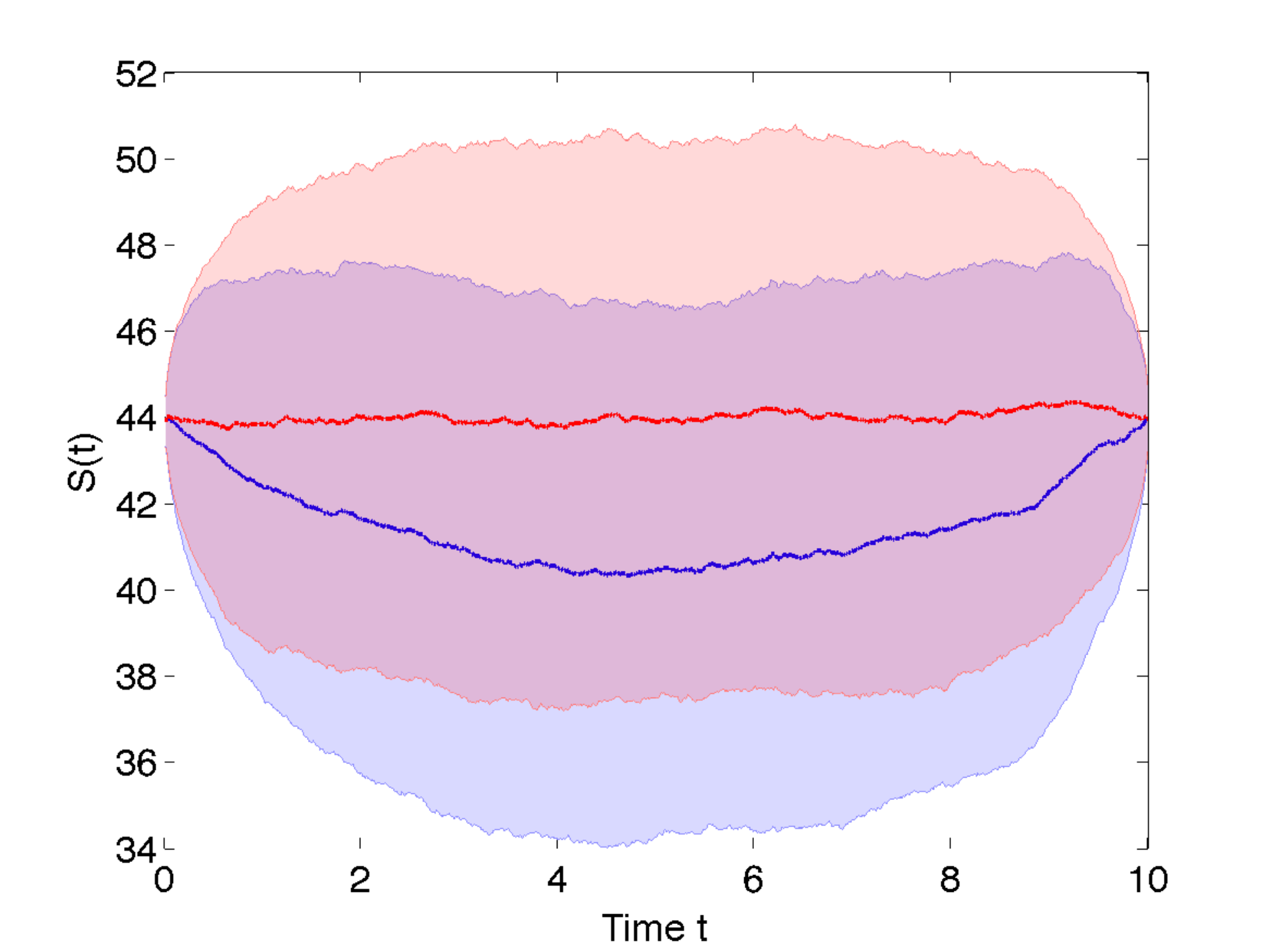}}}}
\caption{\it (a) 10 trajectories of the slow variable conditioned on $S(0) = 44 =
  S(10)$, sampled using the CMA  approximate effective generator. (b) Mean and standard deviation of 1000
  trajectories of the slow variable conditioned on $S(0) = 44 =
  S(10)$, sampled using the approximate effective generator from both the QSSA (blue plots) and CMA (red
  plots).}
\end{center} 
\label{Fig:simple_trajs}
\end{figure} 

Figure \ref{Fig:simple_trajs} (a) shows 10 example trajectories
sampled using the the conditioned path sampling algorithm with the CMA
approximation of the effective generator of the slow variable. We also
implemented exactly the same approach using the QSSA approximation of
the effective generator. The mean and standard deviation of 1000 paths
sampled using both methods is plotted in Figure
\ref{Fig:simple_trajs} (b). Since the paths are conditioned to start and
finish at the mean of the system's monomodal invariant distribution,
we would expect the mean to converge to a constant $S=44$ as we sample
more paths. 

This appears to be the case for the paths sampled using the CMA effective
generator, which is what we would hope since this generator preserves
the true mean of the slow variables, as demonstrated in the previous
section. 

The QSSA, as has also been demonstrated in Section
\ref{sec:lin:qssa}, does not correctly preserve the invariant
distribution of the slow variables, and underestimates the mean value
of the invariant distribution. This can be seen in
\ref{Fig:simple_trajs} (right), where the mean value of the path dips
in the middle of the trajectory as it reverts to the mean of the
invariant distribution of the QSSA approximation, before increasing
towards the end of the trajectory in order to satisfy the condition
$S=44$.

This demonstrates that the accuracy of the approximation of the
effective dynamics has a knock-on impact, as one would expect, to the
accuracy of the conditioned path sampling. It would be preferable,
naturally, if we could compare path statistics of the multiscale
approaches to that of conditioned paths statistics of the full
system. However, this is simply not feasible, due to the size of the
matrices, even for the truncated domain $\Omega$. Instead, this does succeed in demonstrating that these
methods make conditional path sampling of the slow variables a
possibility, where it was computationally intractable previously.  \blue{Since the rates could be explicitly calculated for this
simple example, the effective generators could be produced in the
order of $10^-3$ seconds for the domain $S \in
\{0,1,\ldots,400\}$. The set up process for the path sampling, involving finding the
probabilities in \eqref{eq:S:numR} and computing the required powers of
$\mathcal{M}$ took around 90 seconds. After this, each path took a third of a second to sample.}
}

\red{\subsection{A Bistable Example}\label{sec:bis}
Sampling of paths conditioned on their endpoints is an integral part of
some approaches to Bayesian inversion of biochemical data. A Gibb's
sampler can be used to alternately update the network structure and system parameters, and the
missing data (i.e. the full trajectory), sampled for example using the
method found in \cite{fearnhead2006exact}. However,
efficient methods to sample paths of multiscale systems may also be
useful in other areas. For instance, it may allow us to sample paths
which make rare excursions, or large deviations from mean
behaviour. This forms part of the motivation for considering the next example.

Let us consider the following chemical system, which in certain
parameter regimes exhibits bistable behaviour.
\begin{eqnarray}
R_1, R_2 &:& \qquad\qquad \LRARR{X_2}{X_1+X_2}{k_1}{k_2},
\nonumber \\
R_3, R_4 &:& \qquad\qquad \LRARR{\emptyset}{X_1}{k_3}{k_4},
\label{eq:bis}\\
R_5, R_6 &:& \qquad\qquad \LRARR{X_1 + X_1}{X_2}{k_5}{k_6},
\nonumber \\
R_7 &:& \qquad\qquad \RARR{X_2}{\emptyset}{k_7}.\nonumber
\end{eqnarray}
In particular, we consider parameter regimes where the occurrence of
reactions $R_5$ and $R_6$ are on a relatively faster timescale than
the other reactions. The following is just such a parameter regime:
\begin{gather}
k_1 = 142, \qquad \frac{k_2}{V} = 1, \qquad k_3V = 880,\\\label{eq:bis:params}
k_4 = 92.8, \qquad \frac{k_5}{V} = 10,  \qquad k_6 = 500, \qquad k_7 = 6.\nonumber
\end{gather}
\blue{As before, $V$ denotes the volume of the well-mixed thermally-equilibrated reactor.}

\begin{figure}[htp] 
\begin{center} 
\mbox{
\subfigure[]{\scalebox{0.3}{\includegraphics{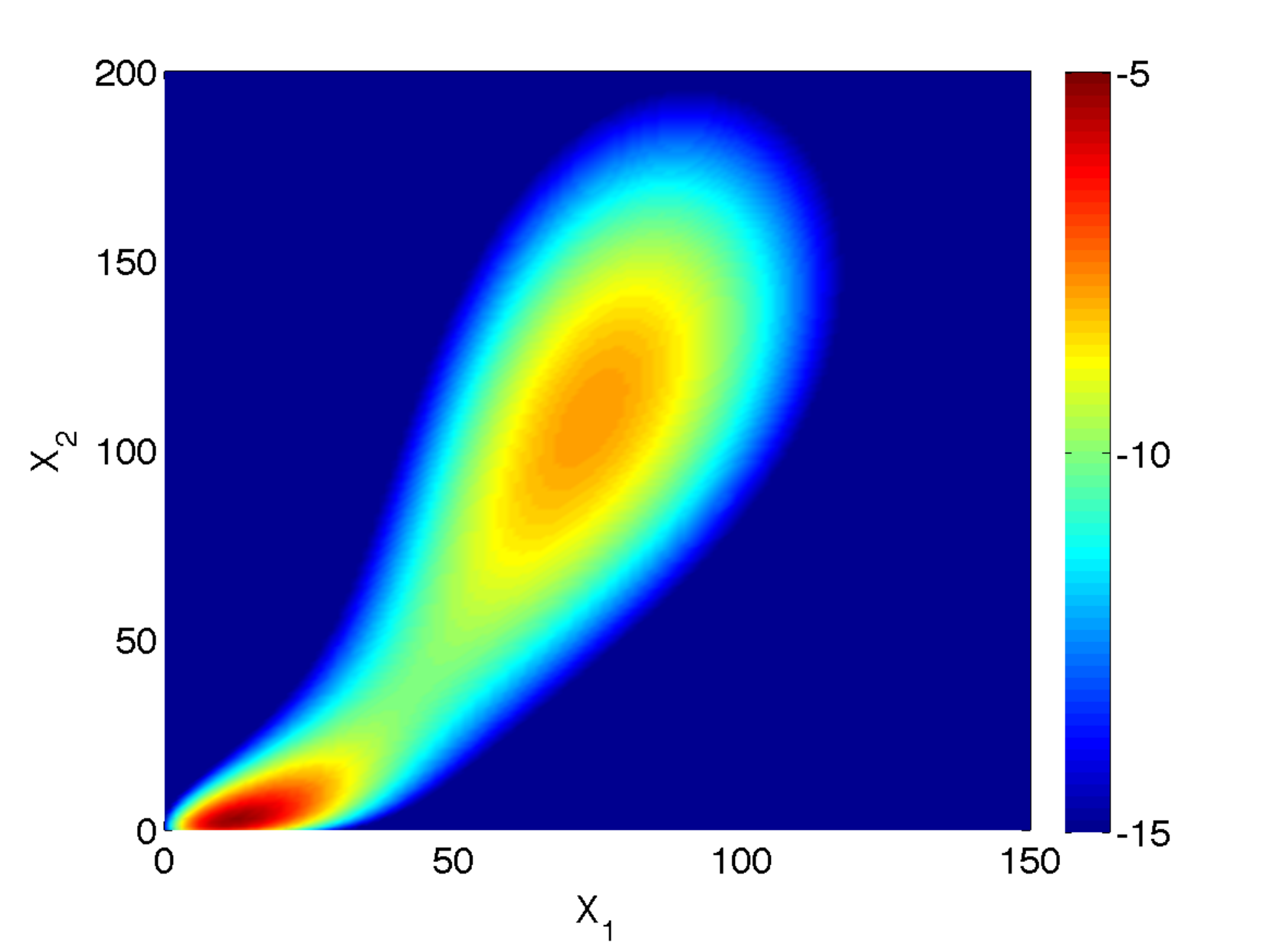}}}
\subfigure[]{\scalebox{0.3}{\includegraphics{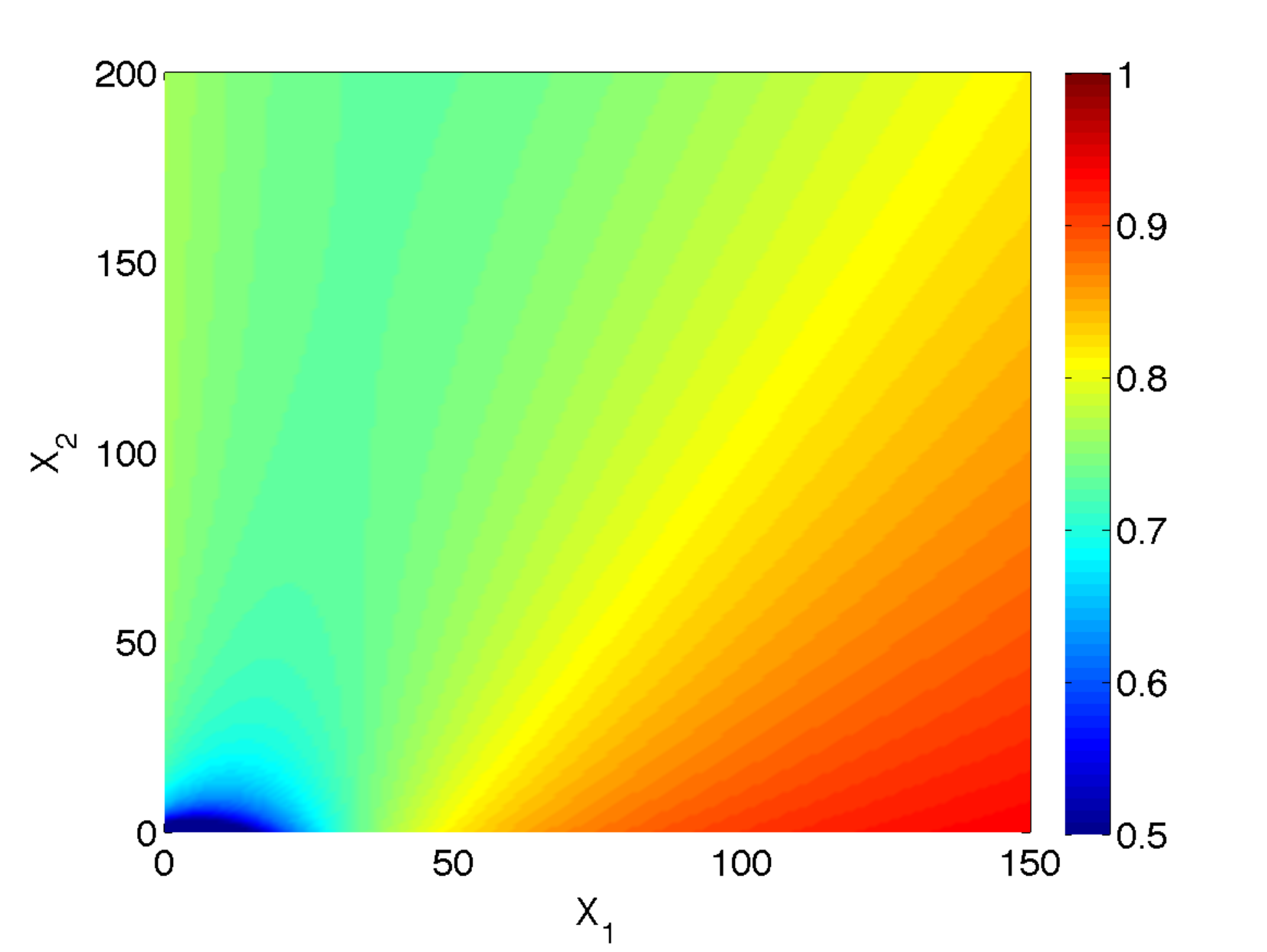}}}}
\caption{\it \blue{(a) A log plot of an approximation $\pi_\Omega$ of the invariant distribution on the
  slow variable $S = X_1 + 2X_2$ of
  system \eqref{eq:bis} with parameters \eqref{eq:bis:params},
  demonstrating the bistable nature of the system. Approximation was
  computed by finding the null space of the full generator of the
  system on the truncated domain
  $\{0,1,\ldots,800\}\times\{0,1,\ldots,1200\}$. (b) Proportion of total propensity $P_{R_5,R_6}(X_1,X_2)$  attributed to the fast
  reactions $R_5$ and $R_6$, given by \eqref{eq:FPP}.\label{fig:bis}}}
\end{center} 
\end{figure} 

That said, this parameter regime is one in which the use of the QSSA will create
significant errors, since the timescale gap is not very large in all
parts of the domain as demonstrated in Figure \ref{fig:bis}. Figure
\ref{fig:bis} (a) shows a highly accurate approximation of the
invariant distribution of the full system, found by computing the
null space of the full
generator for the system truncated to the finite domain $\Omega = \{(x_1,x_2)
\in \{0,1,\ldots,800\}\times\{0,1,\ldots,1200\}$. The
zero eigenvector of this truncated generator was found using
standard eigenproblem solvers, then normalised. Since this system has 2nd order reactions, its invariant
density cannot in general be written in closed form, and as such, we
could use this approximation on the truncated domain in order to quantify the accuracy of the
multiscale approaches. This plot demonstrates the bistable nature of this system,
which can take a long time to switch between the two favourable
regions. This example has
been chosen in order that such an approximation can still be computed
in order to check the accuracy of the approach. 

Figure \ref{fig:bis} (b)
shows the proportion of the total propensity for each state which is
attributed to the fast reactions, $R_5$ and $R_5$, given by:
\begin{equation}\label{eq:FPP}
P_{R_5,R_6}(X_1,X_2) = \frac{\alpha_5(X_1,X_2) + \alpha_6(X_1,X_2)}{\alpha_0(X_1,X_2)} = 
\frac{\alpha_5(X_1,X_2) + \alpha_6(X_1,X_2)}{\sum_{i=1}^{M}\alpha_i(X_1,X_2)}.
\end{equation}
This proportion, which is a measure of the gap in timescales between
the 
``fast'' reactions $R_5$ and $R_6$, and the rest of the reactions,
varies across the domain. We can approximate the expected
proportion of propensity attributed to the fast reactions:
\begin{equation*}
\mathbb{E}(P_{R_5,R_6}) = \sum_{(X_1,X_2) \in \Omega}
P_{R_5,R_6}(X_1,X_2) \pi_\Omega(X_1,X_2),
\end{equation*}
 where $\pi_\Omega$ is the approximate invariant density of the full
 generator on the truncated domain $\Omega$. In this system with
 parameters \eqref{eq:bis:params}, $\mathbb{E}(P_{R_5,R_6}) = 0.6941$,
 i.e. the expected proportion of all reactions which are either of
 type $R_5$ or $R_6$ is $69.41\%$. As such,
 although reactions $R_5$ and $R_6$ are occurring more frequently than
 other reactions, there is not a stark difference in timescales, as we
 might expect in a system for which the QSSA yields a good approximation.
The ``fast'' reactions in this example are reactions $R_5$ and $R_6$, and
as such, $S=X_1 + 2X_2$ is a good choice of slow variable, since this
quantity is invariant to these fast reactions. 

\subsubsection{The QSSA Approach}
By applying the QSSA to the system \eqref{eq:bis}, we can approximate
the effective rates of the slow variables by considering the fast
reactions in isolation. The fast subsystem is given by
the reactions $R_5$ and $R_6$:
\begin{eqnarray}\label{eq:FSS}
C_1 &:& \qquad \qquad X_1 + 2X_2 = S,\\ \nonumber
R_5, R_6&:& \qquad \qquad \LRARR{X_1 + X_1}{X_2}{k_5}{k_6}.
\end{eqnarray}
Lines denoted by $C_i$ in this and what follows denotes a
constraint. It is important to keep a track of these constraints,
since each one reduces the dimension of the state space by one.

For a fixed value of $S=X_1 + 2X_2 \in \{0,1,\ldots,S_{\rm max}\}$, we
wish to find the generator for the process $X_2$ (or equivalently
$X_1 = S - 2X_2$) within this fast subsystem. The generator can be found by considering the master equation
for each state $X_2 = i$:
\begin{eqnarray*}
\frac{dp_i}{dt} &=& -(\alpha_5(S-2i,i) + \alpha_6(S-2i,i))p_{i} +
\alpha_5(S-2(i-1),i-1)p_{i-1} \nonumber \\ &+& \alpha_6(S-2(i+1),i)p_{i+1},
\end{eqnarray*}
where $p_i(t)$ is the probability of $X_2(t) = i$.
Putting this set of differential equations into vector form
gives us:
\begin{equation*}
\frac{d \bP}{dt} = \mathcal{G} \bP,
\end{equation*}
where $\mathcal{G}$ is the generator of the fast subsystem
\eqref{eq:FSS}. Note that since we are restricted to states such that
$X_1 + X_2 = S$ for some value of $S$, there are only $\left \lfloor
  \frac{S}{2} \right \rfloor$ possible different states, and as such
$\mathcal{G} \in \mathbb{R}^{\left \lfloor
  \frac{S}{2} \right \rfloor \times \left \lfloor
  \frac{S}{2} \right \rfloor}$. Even for moderately large values of
$S$, the one-dimensional null space of such a sparse matrix is not computationally
expensive to find, and when normalised gives us the invariant density of $X_2$ (and
therefore $X_1$ if required). This invariant density can then be used
to compute the expectation of the propensities of the slow reactions
of the system for the state $S$ as in \eqref{effprop}, and in turn be entered into the
(truncated) effective generator for the slow variable.

\subsubsection{The Constrained Approach}
When using the CMA, the methodology is largely the same as was
described for the QSSA-based approach in the last section. The only real
difference lies in the subsystem which is analysed in order to compute
the invariant distribution of the fast variables conditioned on the
value of the slow variable. As we have done previously, we will
consider each of the reactions in the system in turn, constraining the
value of the slow variable to a particular value, whilst being sure
not to change the value of the fast variables. There are two choices
for the fast variable, in order to form a basis of the state space
along with the slow variable $S$, but as explained in detail in
\cite{cotter2011constrained}, $F = X_2$ is the best choice, since
there is a zeroth order reaction involving $X_1$, which can lead to
an unphysical constrained subsystem, if this is chosen as the fast variable.

With this choice of fast variables, the first four reactions all
disappear in the constrained subsystem. This is because none of these
reactions alter the fast variable, and as such the constrained
stoichiometric projector maps their stoichiometric vectors to zero,
and therefore reactions $R_1, R_2, R_3, R_4$ have no net effect on
the constrained subsystem.

Reaction $R_7$ differs in that it causes a change in the fast variable
$X_2$. The projector in this case maps the stoichiometric vector to
$[-2,1]^T$ and therefore the net effect of reaction $R_7$ is
equivalent to $\RARR{X_2}{X_1 + X_1}{k_7}$. This leads to the
following constrained system:
\begin{eqnarray*}
C_1&:& \qquad \qquad X_1 + 2X_2 = S,\nonumber\\
R_5, R_6&:& \qquad\qquad \LRARR{X_1 + X_1}{X_2}{k_5}{k_6},\\
R_7&:& \qquad\qquad \RARR{X_2}{X_1 + X_1}{k_7}.\nonumber
\end{eqnarray*}
Note that since reactions $R_6$ and $R_7$ have the same stoichiometry,
this system can be simplified by removing $R_7$ and adding its rate to
$R_6$:
\begin{eqnarray}\label{eq:bis:con}
C_1&:& \qquad \qquad X_1 + 2X_2 = S,\nonumber\\
R_5, R_6&:& \qquad\qquad \LRARR{X_1 + X_1}{X_2}{k_5}{k_6 + k_7}.
\end{eqnarray}
For every fixed value of $S \in \{0,1,\ldots,S_{\rm max}\}$, the
generator for \eqref{eq:bis:con} can be found following the same
approach as in the previous section, the only difference being the
altered rate for reaction $R_6$. Following this methodology, an effective generator $\mathcal{G}$ can
be computed.

\subsubsection{Comparison of approximation of invariant densities}
\begin{figure}[htp] 
\begin{center} 
\includegraphics[width=0.6\textwidth]{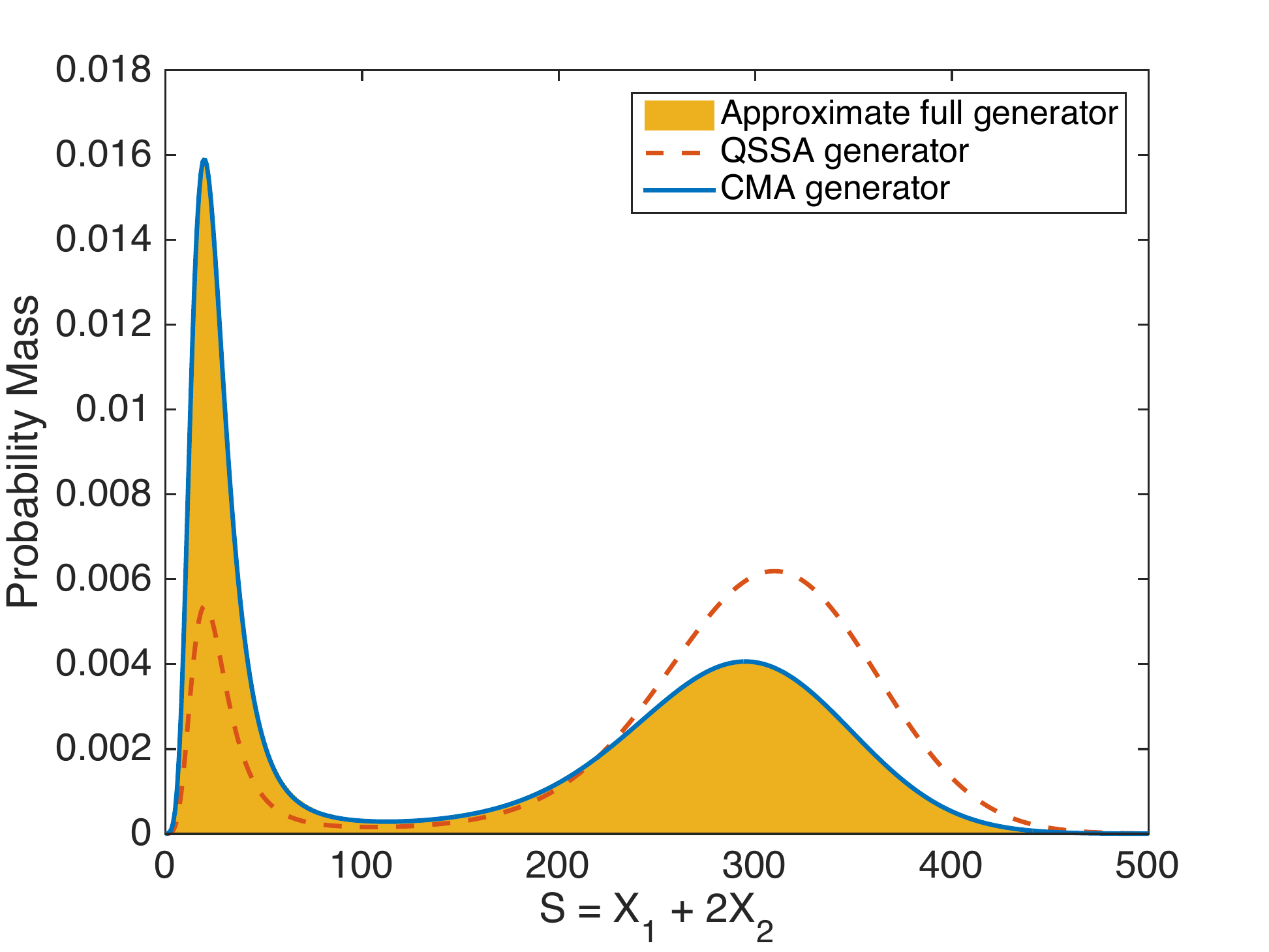}
\caption{\it Approximations of the invariant
  distribution of the slow variable $S=X_1 + 2X_2$ of system
  \eqref{eq:bis} with parameters \eqref{eq:bis:params}, through
  computing the null space of the truncated generator of the full
  system (histogram), of the effective generator computed using the
  CMA (solid line) and computed using the QSSA (dashed line).\label{Fig:bis2}}
\end{center} 
\end{figure} 

One approach to quantifying the accuracy of these two methods of
approximating effective generators of the slow variable, is to compare
the invariant distributions of the two systems with that of the
marginalised density of the slow variable in the full system. We
consider the approximation $\pi_\Omega$ of the invariant density of the full
system, truncated to the region $\Omega = \{(x_1,x_2)
\in \{0,1,\ldots,800\}\times\{0,1,\ldots,1200\}$, as shown in Figure
\ref{fig:bis} (a). We can marginalise this density to
find an approximation of the invariant density of the slow variable,
as is shown by the histogram in Figure \ref{Fig:bis2}.

\begin{table}[]
\centering
\begin{tabular}{l|ccc|}
\cline{2-4}
                                                & QSSA                   & CMA                   & $\pi_\Omega$          \\ \hline
\multicolumn{1}{|l|}{Relative $l^2$ difference} & $6.347\times 10^{-1}$     & $1.796\times 10^{-2}$ & -                     \\
\multicolumn{1}{|l|}{LH peak position}          & $20$                   & $20$                  & $20$                  \\
\multicolumn{1}{|l|}{LH peak height}            & $5.378\times 10^{-3}$  & $1.591\times 10^{-2}$ & $1.582\times 10^{-2}$ \\
\multicolumn{1}{|l|}{RH peak position}          & $309$                  & $295$                 & $295$                 \\
\multicolumn{1}{|l|}{RH peak height}            & $6.192 \times 10^{-2}$ & $4.060\times 10^{-3}$ & $4.006\times 10^{-3}$ \\ \hline
\end{tabular}
\caption{{\it Differences in the accuracy of the QSSA and CMA approximations of the
    invariant density of $S$, with respect to the approximation $\pi_\Omega$.}}
\label{tab:bis}
\end{table}

The CMA approximation of the invariant density of the slow variable is
indistinguishable by eye from the highly accurate approximation
computed in this manner, as shown in Figure \ref{Fig:bis2}. The QSSA
approximation, on the other hand, incorrectly approximates both the
placement and balance of probability mass of the two peaks in the
distribution. The difference in the quality of these approximations is
stark. This example is an extreme one, as the parameters have been
chosen to demonstrate how far apart these two approximations can be,
but since the CMA has no additional costs associated with it, the
advantages of this approach are significant. The relative $l^2$ errors
of these two approaches, when compared with the approximate density
$\pi_\Omega$, are given in Table \ref{tab:bis}, along with the
position and heights of the two local maxima in the densities.

\blue{The CMA computed the generator on the domain $S\in[0,2000]$ in around
55 seconds, and the eigensolver took less than a tenth of a second to
find the null space to approximate the invariant density. This is
negligible in comparison with the cost of exhaustive stochastic
simulation of the full system.}

\subsubsection{Conditioned path sampling using effective generators}
Given an approximation of the effective generator of the slow
variables, computed using the CMA or the QSSA, we can now employ the methodology
of \cite{fearnhead2006exact}, as summarised in Section \ref{sec:CPS},
to sample paths conditioned on their endpoints. This time, a full
eigenvalue decomposition of the matrix $\mathcal{M} =
\frac{1}{\rho}\mathcal{G} + I$ was computed, so that matrices $V$ and
$D$ could be found with $V$ unitary and $D$ diagonal, with
$\mathcal{M} = V^{-1}DV$. Then rows of $\mathcal{M}^r = V^{-1}D^rV$ can be
efficiently and accurately computed, as required in \eqref{eq:S:numR}
and \eqref{eq:S:type}.

\begin{figure}[htp] 
\begin{center} 
\mbox{ 
\subfigure[]{\scalebox{0.3}{\includegraphics{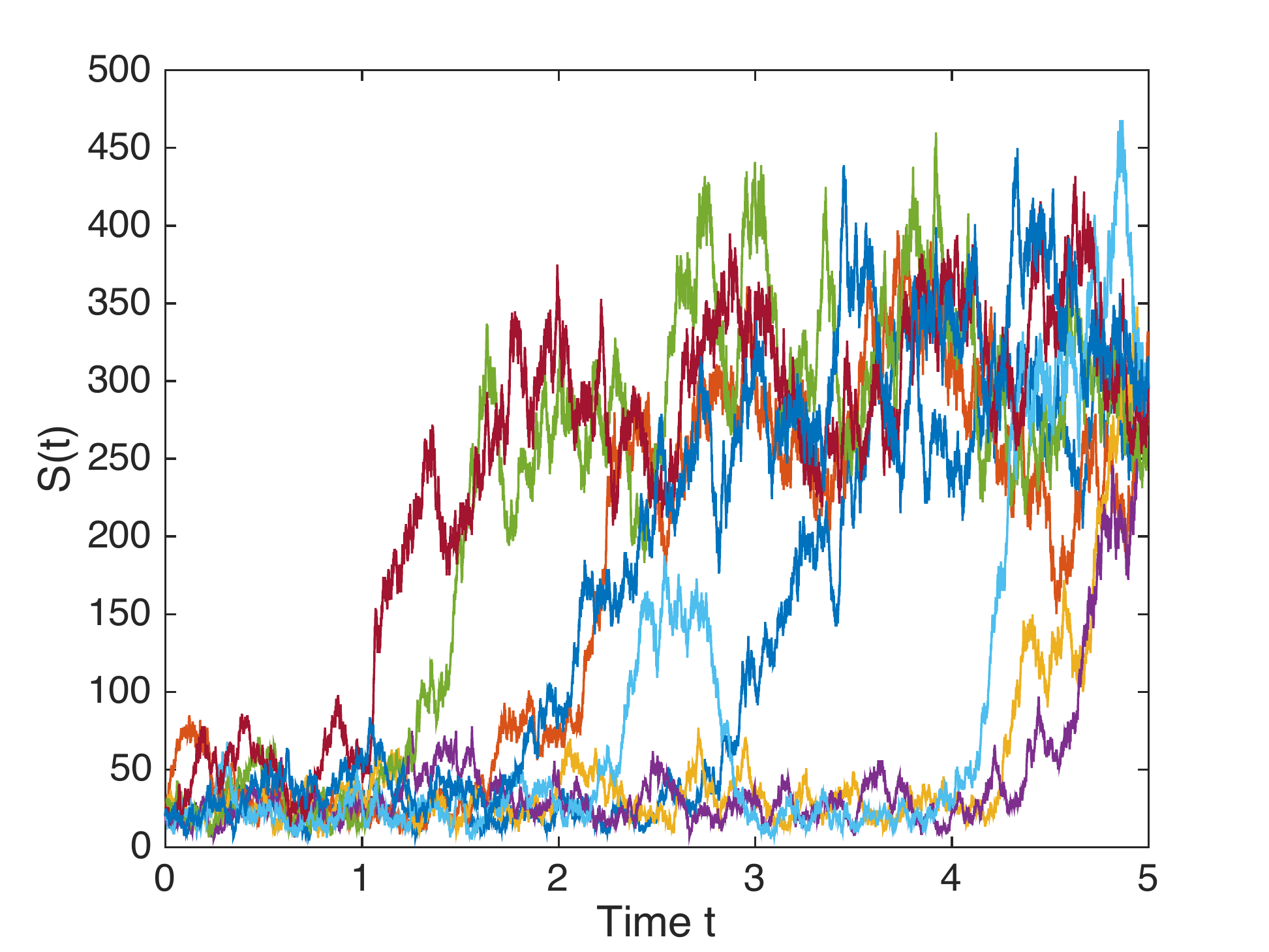}}} 
\subfigure[]{\scalebox{0.3}{\includegraphics{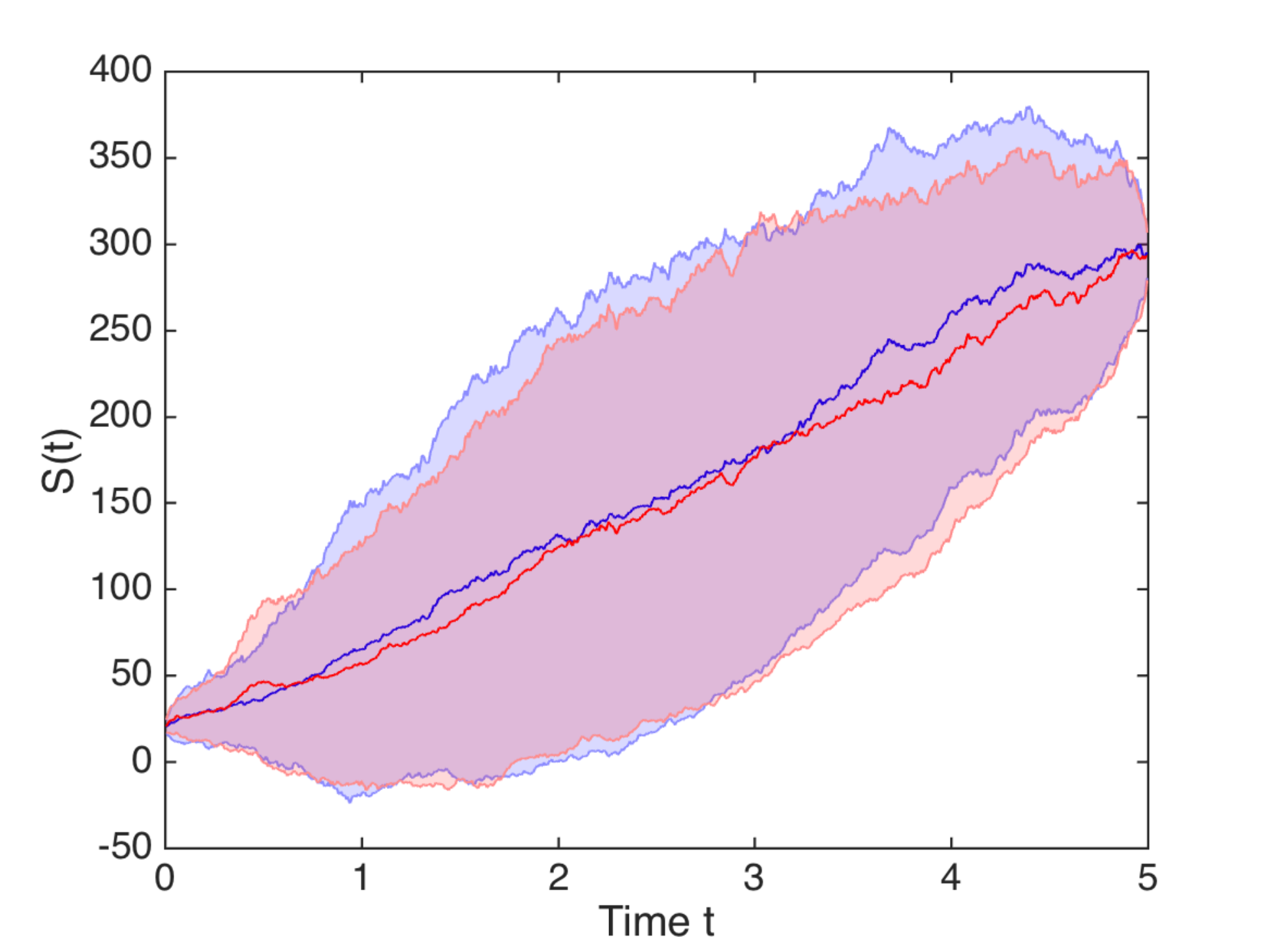}}} 
}
\caption{\it (a) 8 trajectories of the slow variable $S=X_1+2X_2$
  sampled conditioned on $S(0) =
20,\, S(10) = 195,\, S(t) \in \Omega = \{0,1,\ldots,500\} \forall t\in [0,5]$ for the
system \eqref{eq:bis} with parameters \eqref{eq:bis:params}, using the
CMA approximation of the effective generator. (b) The means and
standard deviations of 100 paths sampled using the QSSA (blue plots)
and CMA (red plots).\label{Fig:bis_trajs}}
\end{center} 
\end{figure}

Figure \ref{Fig:bis_trajs} presents results using this approach. An
effective generator for the system \eqref{eq:bis} was computed for the
domain $X_1 + 2X_2 = S \in \Omega = \{0,1,\ldots,500\}$, using both
the QSSA and CMA,  and then fed
into the conditioned path sampling algorithm. Figure
\ref{Fig:bis_trajs} (a) shows 8 samples of conditioned
paths approximated using the CMA. Notice that as the transition time between the two favourable
regions is relatively short compared with the length of the simulation,
the time of the transition varies greatly between the different
trajectories. This indicates that we are producing trajectories with a
fair reflection of what happens in a transition between these
regions. Figure \ref{Fig:bis_trajs} (b) shows the means and
standard deviations of 100 paths sampled for both methods of computing
the effective generator. The QSSA, which overestimates the value of
the second peak in the invariant density, has a
higher mean than the CMA. This demonstrates again that errors in
approximating the effective generator has a knock-on affect to
applications such as conditioned path sampling.}

\blue{The effective generator was computed on the domain $S\in
  [0,500]$ for the path sampling, which took the CMA close to 5
  seconds to approximate. The calculation of the probabilities in
  \eqref{eq:S:numR}, and the full eigenvalue decomposition of the
  generator matrix on this domain, took around 50 seconds. After this,
  each path took around 350 seconds to sample.}

\subsection{An Example of the Nested CMA Approach}\label{sec:num:nest}

In this section, we will illustrate how the nested approach outlined
in Section \ref{sec:Nest} can be applied. \red{We will consider an example
for which we know the invariant distribution of the slow variables. This gives us a way of quantifying any errors that we incur
by applying the nested CMA and QSSA approaches.}
\begin{eqnarray}
R_1 &:& \qquad\qquad \emptyset \, \overset{k_1}{\longrightarrow} \, X_1,
\nonumber \\
R_2 &:& \qquad\qquad X_3 \, \overset{k_2}{\longrightarrow} \, \emptyset,
\nonumber \\
R_3 &:& \qquad\qquad X_1 \, \overset{\kappa}{\longrightarrow} \, X_2,
\label{eq:lin2} \\
R_4 &:& \qquad\qquad X_2 \, \overset{\kappa}{\longrightarrow} \, X_1,
\nonumber  \\
R_5 &:& \qquad\qquad X_2 \, \overset{\gamma}{\longrightarrow} \, X_3,
\nonumber \\
R_6 &:& \qquad\qquad X_3 \, \overset{\gamma}{\longrightarrow} \, X_2.  
\nonumber 
\end{eqnarray}
We will consider this system in the following parameter regime:
\red{\begin{equation} \label{eq:lin2:params}
k_1V = 20, \qquad k_2 = 1, \qquad \kappa = 100, \qquad \gamma = 10.
\end{equation}} 
\blue{As before, $V$ denotes the volume of the well-mixed thermally-equilibrated reactor.}
In this regime, there are multiple different time scales on which the
reactions are occurring. This is demonstrated in Figure
\ref{fig:lin2:occ}, where there is a clear gap in the frequency of
reactions $R_1$ and $R_2$ (the slowest), $R_5$ and $R_6$ (fast
reactions) and $R_3$ and $R_4$ (fastest reactions).

\red{\begin{figure}[htp] 
\begin{center} 
\includegraphics[width=0.6\textwidth]{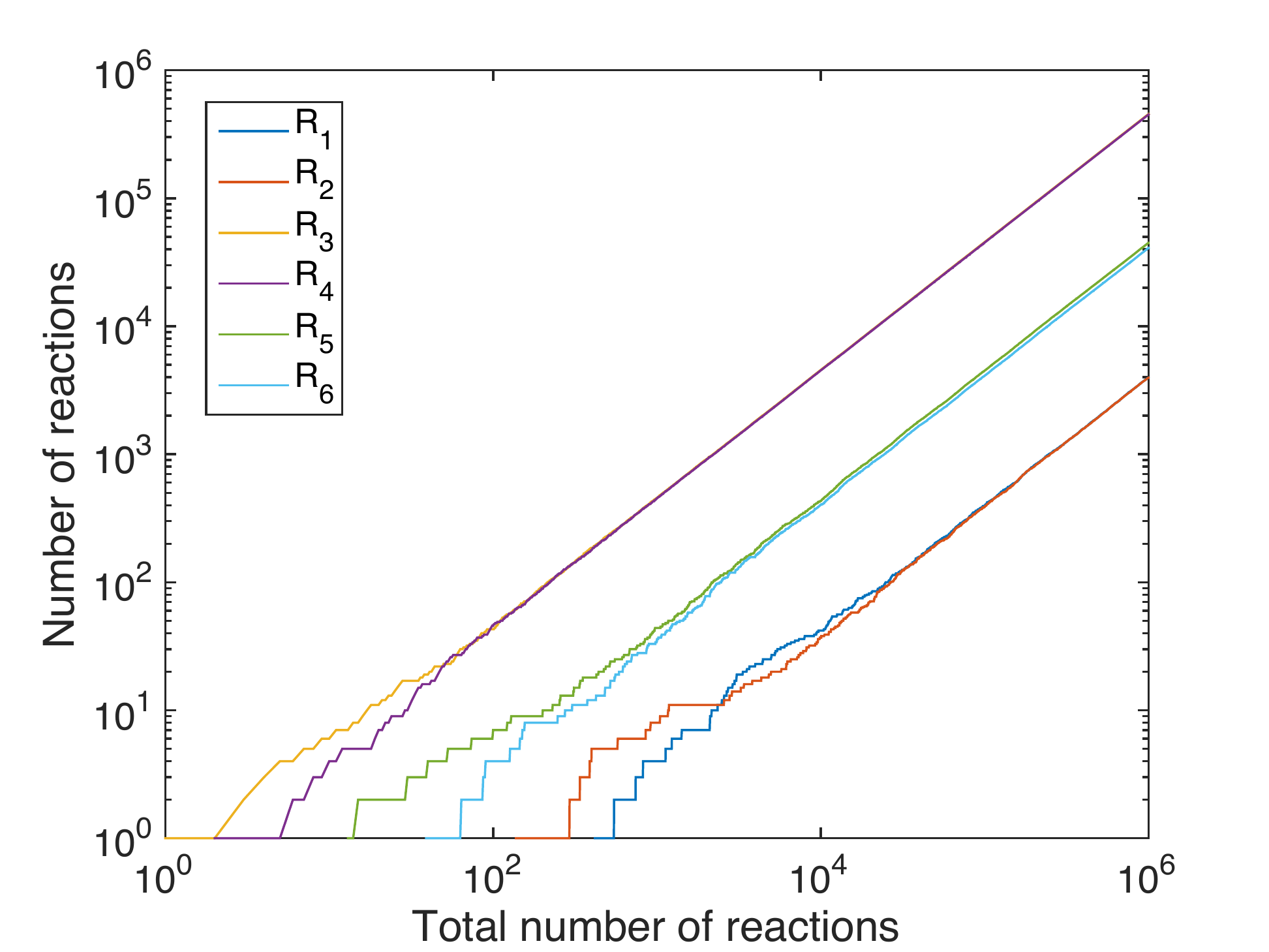}
\caption{\it Relative occurrences of the reactions $R_1$-$R_6$, for
  the system \eqref{eq:lin2} with parameters
  \eqref{eq:lin2:params}. The most frequent reactions are reactions
  $R_3$ and $R_4$, reactions $R_5$ and $R_6$ are the next most
  frequent, with reactions $R_1$ and $R_2$ being the least frequent. \label{fig:lin2:occ}}
\end{center} 
\end{figure} }
\red{
This system was
chosen as we are able to, using the results in
\cite{jahnke2007solving}, find the exact invariant distribution of
the slow variable $S_1 = X_1 + X_2 + X_3$. In this instance, it is a Poisson
distribution with mean parameter $$\lambda =
\frac{k_1V}{k_2\gamma\kappa} \left ( \gamma k_2 + 3\gamma \kappa +
  2k_2\kappa \right ) =
64.2.$$

\subsubsection{QSSA-based analysis}
One method to analyse such a system would be a nested QSSA-based
analysis, similar to that which is suggested in
\cite{weinan2005nested}. In this paper the authors consider systems
with reactions occurring on multiple timescales. If at first we
consider  all reactions $R_3$-$R_6$ to be fast reactions, then by
applying the QSSA we are interested in finding the invariant
distribution of the following
fast subsystem:

\begin{eqnarray}
C_1 &:& \qquad\qquad X_1 + X_2 + X_3 = S_1,\nonumber\\
R_3 &:& \qquad\qquad X_1 \, \overset{\kappa}{\longrightarrow} \, X_2 , \nonumber
\\
R_4 &:& \qquad\qquad X_2 \, \overset{\kappa}{\longrightarrow} \, X_1,
\label{eq:lin2:f1}  \\
R_5 &:& \qquad\qquad X_2 \, \overset{\gamma}{\longrightarrow} \, X_3,
\nonumber \\
R_6 &:& \qquad\qquad X_3 \, \overset{\gamma}{\longrightarrow} \, X_2.  
\nonumber 
\end{eqnarray}

Note that the quantity $S_1 = X_1 + X_2 + X_3$ is a conserved quantity
with respect to these reactions, and as such is the slow variable in
this system. This is in itself also a system with more than one timescale, and as
such, we may want to iterate again and apply a second QSSA assumption,
based on the fact that reactions $R_3$ and $R_4$ are fast reactions in
comparison with reactions $R_5$ and $R_6$. This leads to a second fast
subsystem:

\begin{eqnarray*}
C_1 &:& \qquad\qquad X_1 + X_2 + X_3 = S_1,\\
C_2 &:& \qquad\qquad X_1 + X_2 = S_2,\\
R_3 &:& \qquad\qquad X_1 \, \overset{\kappa}{\longrightarrow} \, X_2,
\\
R_4 &:& \qquad\qquad X_2 \, \overset{\kappa}{\longrightarrow} \, X_1.
\end{eqnarray*}

Note that the quantity $S_2 = X_1+X_2 = S_1 - X_3$ is  a conserved quantity
with respect to these reactions, and as such is the slow variable in
this system. At this point in \cite{weinan2005nested}, the authors
simulate the system using the Gillespie SSA. We could adopt the
approach that we used in Section \ref{sec:bis}, in which we find the
invariant distribution of the system by constructing its generator and
then finding the normalised eigenvector corresponding to the null space
of that operator. This would not be expensive since there are only
$S_2$ different states. However, as in Section \ref{sec:lin}, as this
system only contains monomolecular reactions, we can exactly find its
invariant distribution. In this case, $X_1$ and $X_2$ follow a  binomial
distribution with mean $\bar{X_1},\bar{X_2} = \frac{S_2}{2}$. This can
then be used to compute the effective rate of reaction $R_5$ in the
first subsystem \eqref{eq:lin2:f1}, $\alpha_5(X_1,X_2) \approx \gamma\bar{X_2} =
\frac{\gamma}{2}S$. This fast subsystem is then reduced to the following:

\begin{eqnarray*}
C_1 &:& \qquad \qquad X_1 + X_2 + X_3 = S_1 = S_2 + X_3,\\
C_2 &:& \qquad \qquad X_1 + X_2 = S_2,\\
R_5 &:& \qquad\qquad S_2 \, \overset{\gamma/2}{\longrightarrow} \, X_3,
\nonumber \\
R_6 &:& \qquad\qquad X_3 \, \overset{\gamma}{\longrightarrow} \, S_2.  
\nonumber 
\end{eqnarray*}
Note that we have completely eliminated the fast variables $X_1$ and
$X_2$, and instead consider the slower variable $S_2=X_1 + X_2$, with
effective rate for $R_5$ given by the analysis above. This system is
exactly solvable, and its invariant distribution is a gamma distribution with
means given by $\bar{X_3} = \frac{S_1}{3}$ and $\bar{S_2} =
\frac{2S_1}{3}$, found by computing the steady states of the mean
field ODEs\cite{jahnke2007solving}. This in turn can be used to
compute the effective rate of reaction $R_2$ in the full system, where we
now lose all of the fast variables $X_1,X_2,X_3$ and instead wish to
understand the dynamics of the slow variable $S_1 = X_1 + X_2 + X_3$,
which is only altered by reactions $R_1$ and $R_2$. This system is
given by the following:

\begin{eqnarray*}
R_1 &:&  \qquad\qquad \emptyset \, \overset{k_1}{\longrightarrow} \, S_1,
\nonumber \\
R_2 &:&  \qquad\qquad S_1 \, \overset{k_2/3}{\longrightarrow} \, \emptyset.
\nonumber \\
\end{eqnarray*}
Here the effective rate for $R_2$ has been found by using the
approximation of the effective rate
$\alpha_2(S_1) = k_2\bar{X_3} = \frac{k_2}{3}S$.
\subsubsection{CMA-based analysis}
We will now go through the same procedure, but this time using the
constrained subsystems instead of the fast subsystems as used in the
previous section.} There are 3
choices for the fast reactions, each involving two out of $X_1$, $X_2$
and $X_3$. Since $X_1$ is the product of a zeroth order reaction, it
is preferable not to include this as one of the fast variables, and so
we pick ${\bf F}_1 = [X_2,X_3]$. We then construct the constrained
subsystem for this choice of slow and fast variables:
\begin{eqnarray}
C_1 &:& \qquad\qquad X_1 + X_2 + X_3 = S_1, \nonumber \\
R_2 &:& \qquad\qquad X_3 \, \overset{k_2}{\longrightarrow} \, X_1,
\nonumber \\
R_3 &:& \qquad\qquad X_1 \, \overset{\kappa}{\longrightarrow} \, X_2,
\nonumber \\
R_4 &:& \qquad\qquad X_2 \, \overset{\kappa}{\longrightarrow} \, X_1,  
\label{eq:lin2:con1}  \\
R_5 &:& \qquad\qquad X_2 \, \overset{\gamma}{\longrightarrow} \, X_3,  
\nonumber \\
R_6 &:& \qquad\qquad X_3 \, \overset{\gamma}{\longrightarrow} \, X_2.  
\nonumber 
\end{eqnarray}
Note that $R_1$ is removed, since it does not change the fast
variables. $R_2$ is the only other reaction which has changes to its
stoichiometry due the constrained stoichiometric projector. We have
reduced the dimension of the system (due to the constraint $X_1 + X_2
+ X_3 = \sigma$ for some $\sigma \in \mathbb{N}$), but we are still left with a multiscale system, which in
theory could be computationally intractable for us to find the
invariant distribution for, through funding the null space of its
generator. Therefore, we can apply another iteration of the CMA to
this constrained system.

Reactions $R_3$ and $R_4$ are the fastest reactions in the system, and
therefore we pick our next slow variable that we wish to constrain to
be $S_2 = X_1 + X_2$, which is invariant with respect to these
reactions. Due to the previous constraint $S_1 = X_1 + X_2 + X_3$, we are only required to
define one fast variable for this system. \red{Both choices $F_2 = X_1,
X_2$, are essentially equivalent, and so we pick $F_2 = X_1$}. These
choices lead us to the following second constrained system:
\begin{eqnarray}
C_1 &:& \qquad\qquad X_1 + X_2 + X_3 = S_1, \nonumber \\
C_2 &:& \qquad\qquad X_1 + X_2 = S_2, \nonumber \\
R_2 &:& \qquad\qquad \alpha_2({\bf X}) = \begin{cases} k_2X_3, \qquad
  &{\rm if} \, X_2 > 0, \\ 0 & {\rm
    otherwise,} \label{eq:lin2:con2}\end{cases} \\ & & \qquad \qquad {\boldsymbol{\nu}}_2
= [1,-1,0]^T, \nonumber
\\
R_3 &:& \qquad\qquad X_1 \, \overset{\kappa}{\longrightarrow} \, X_2,  
\nonumber\\
R_4 &:& \qquad\qquad X_2 \, \overset{\kappa}{\longrightarrow} \, X_1.  
\nonumber 
\end{eqnarray}
\blue{Here ${\boldsymbol{\nu}}_i$ denotes the stoichiometric vector associated with
  reaction $R_i$, i.e. the vector which is added to the state ${\bf
    X}(t)$ if reaction $R_i$ fires at time $t$.}
\red{Notice that we now have two separate constraints, and as such
reactions $R_5$ and $R_6$ now have zero stoichiometric
vectors. Moreover, these constraints lead us to a somewhat unphysical
reaction for $R_2$. The reactant for this reaction is $X_3$, but only
$X_2$ and $X_1$ are affected by this altered reaction. In system
\eqref{eq:lin2:con1} when reaction $R_2$
fires, we lose one $X_3$, and gain $X_1$. Therefore, both
constraints within \eqref{eq:lin2:con2} have been violated. In order to reset these constraints,
without changing the fast variable $F = X_3$, we arrive at the
stoichiometry presented in \eqref{eq:lin2:con2}. Note that we add the
condition that this reaction can only happen if $X_2> 0$, as we cannot
have negative numbers of any species.}

This is a closed system, with a very limited number of different
states. Therefore, it is computationally cheap to construct its
generator, and to find that generator's null space. Our aim with this
system, is to find the invariant distribution of the fast variable
given particular values for the constraints $C_1$ and $C_2$. This distribution will
then allow us to compute the expectation of the reaction $R_4$ within
the constrained system \eqref{eq:lin:con}, which
is the only reaction which is dependent on the
results of the second constrained system (since $X_3 = S_1 -
S_2$). Once the invariant distribution has been found, this can be
used to find the effective propensity of reaction $R_5$ given values
of $S_1 = X_1 + X_2 + X_3$ and
$S_2 = X_1 + X_2$. In turn, the constrained system
\eqref{eq:lin2:con1} can then be solved to find the invariant distribution
on $X_3$ given a value of $S_1$. Finally, this leads us to the
construction of an effective generator for the slow variable
$S_1$.
 
\red{
Since this final constrained subsystem is aphysical, we cannot use the
results of \cite{jahnke2007solving} to find the invariant
distribution, and as such we must approximate them through finding the
null space of the generator, as we did in Section \ref{sec:bis}

\subsubsection{Comparison of approximation of invariant densities}

\begin{figure}[htp] 
\begin{center} 
\includegraphics[width=0.6\textwidth]{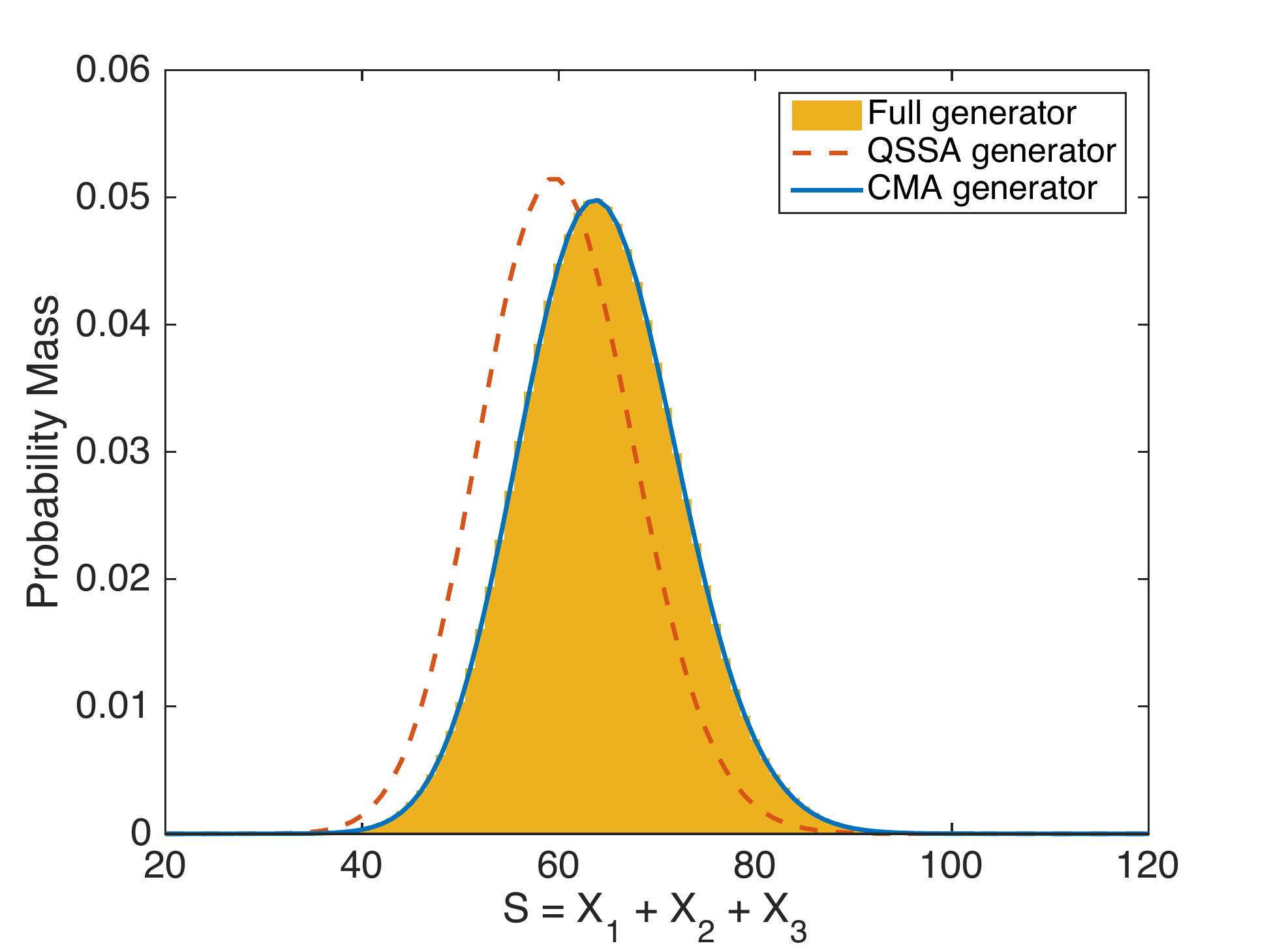}
\caption{\it Approximations of the invariant
  distribution of the slow variable $S=X_1 + X_2 + X_3$ of system
  \eqref{eq:lin2} with parameters \eqref{eq:lin2:params}, through
  marginalisation of the invariant distribution of the full
  system (histogram), of the effective generator computed using the
  CMA (solid line) and computed using the QSSA (dashed line).\label{Fig:lin2}}
\end{center} 
\end{figure} 

Figure \ref{Fig:lin2} shows the invariant distributions of the
slow variables $S=X_1 + X_2 + X_3$ computed by marginalising the invariant
distribution of the full system, and from the CMA and QSSA as outlined
above. The distribution computed using the CMA is indistinguishable by
eye from the true distribution, and has a relative error of
$5.936\times 10^{-12}$, which can be largely attributed to rounding
errors and error tolerances in the eigenproblem solvers. The QSSA approximation, on the other hand, has a significant
relative error of $3.739\times 10^{-1}$. This demonstrates again the
substantial improvements in accuracy that we gain in using the
constrained approach rather than one based on the QSSA. This is
delivered at no substantial additional computational effort. As in the
previous two examples, the highly accurate effective generator
approximated using the CMA can be used in a host of applications where
the full generator could not, such as conditioned path sampling.
}

\blue{The CMA is more expensive in this example than the previous ones, as
  there are a very large number of small eigenvalue problems to
  solve. This is due to the fact that there are reactions of three species
  occurring on three different time scales. The generation of the CMA
  approximation of the effective generator took around 1240 seconds,
  and the subsequent approximation of the invariant distribution of
  the slow variables took just over half a second. This still pales
  into comparison with the cost of exhaustive stochastic simulation of
  the system. The savings would be even more pronounced in systems
  with multimodal invariant distributions where switches between the
  modes are rare.}

\section{Conclusions}\label{sec:conc}
In this paper, we presented a significant improvement and extension to the original
constrained multiscale algorithm (CMA). Through constructing and
finding the null space of the generator of the constrained process, we
can find its invariant distribution without the need for expensive
stochastic simulations. The CMA in this format can also be used not
just to
approximate the parameters of an approximate diffusion, but to
approximate the rates in an effective generator for the slow
variables.

\red{In this paper we have not discussed how the slow and fast variables
  in these systems can be identified. In the simple examples
  presented, this is relatively straightforward. However in
  general, this is far from the case. If the high probability regions
  in the statespace are known a priori, or possibly identified through
  short simulations of the full system, then it is possible to
  identify which are the fast reactions in the system, and therefore
  what good candidates for the slow variable(s) could be. Other more
  sophisticated approaches exist, for example methods for automated
  analysis to identify the slow
  manifold\cite{erban2007variable,singer2009detecting,sarich2012approximating}. One
  relatively ad hoc approach might be to briefly simulate the full
  system using the Gillespie SSA, which can give a good indication as
  to which the fast reactions are. Good candidates for slow variables
  are often linear combinations of the species who are invariant to
  the stoichiometry of the fast reactions, as we have seen in this
  paper. If the regions which the system is highly likely to spend the
  majority of its time are known, then looking at the relative values
  of the propensity functions, as we did in Figure \ref{fig:bis} (b),
  can lead to an understanding of which reactions are fast and which
  are slow.
}

Through iterative nesting, the CMA can be applied to much more complex
systems, as it can be applied repeatedly if the resulting constrained
system is itself multiscale. This makes it a viable approach
for a bigger family of (possibly biologically relevant) systems. This
nested approach breaks up the original task of solving an eigenvalue
problem for one large matrix per row of the effective generator, down
into many eigenvalue solves for significantly smaller generators for
smaller dimensional problems,
making the overall problem computationally tractable.

\red{
In the first example, we demonstrated that the CMA produces an
approximation of the dynamics of the marginalised slow process in the
system which is exact, at least by the measures that we have applied
thus far, in the case of systems of monomolecular reactions. Since
such systems are well understood, we were also able to compare this
with the accuracy of the equivalent QSSA-based method, which incurred
significant errors. We then applied the method of Fearnhead and
Sherlock\cite{fearnhead2006exact} to the approximate effective
generators of the two approaches, in order to approximately sample
conditioned paths of the slow variables. This task would be
computationally intractable to attempt with the full generator for
this system. This also demonstrated how the accuracies of the two
approximations can impact the accuracy of any application for which
they may be used.

In the second example, a more complex bistable system was also analysed using the CMA, and
the invariant distribution of the computed effective generator was
shown to be very close to the best approximation that we could make of
the invariant distribution of the slow variables, using the null space
of the original generator with as little truncation as we could
sensibly manage with our computational resources. This was in stark
contrast with the poor approximation which was computed using the
equivalent QSSA-based approach. This highlighted again the
improvement, at no or little extra cost, of using the constrained
approach as opposed to the QSSA. 

In the final example, we demonstrated how the constrained approach
might be applied to a more complex example with multiple
timescales. The algorithm can be applied iteratively in order to
reduce the constrained subsystems themselves into a collection of
easily solved one-dimensional problems. When comparing the invariant
distributions of the approximate processes computed using the two approaches, the QSSA once again was incorrectly approximating the
distribution of the fast variables conditioned on the slow variables,
and so incurred significant errors. In contrast, the CMA produced an
approximation to the invariant measure which was accurate up to
12 digits.
}

We showed how these effective generators can be used in the sampling
of paths conditioned on their endpoints. Such an approach could be
employed as a method to sample missing data within a Gibb's sampler
when attempting to find the structure of a network that was
observed\cite{fearnhead2006exact}. This approach could also be used
simply to simulate trajectories of the slow variables, in the same
vein as \cite{cao2005slow} or \cite{weinan2005nested}. In this
instance, it would only be necessary to compute the column of the
effective generator corresponding to the current value of the slow
variables.

\red{
The constrained approach consistently significantly outperforms
approximations computed using the more standard QSSA-based approach,
and at negligible additional cost. Furthermore, in the limit of large
separation of timescales, the constrained approach asymptotically
approaches the QSSA approximation.}

\blue{The computational savings that we make in using the CMA depends
  on the application with which we
wish to use the effective generators. Similarly, if we wish to approximate the invariant distribution of the
slow variables, then the CMA will always be less costly than exhaustive stochastic simulation.
This is because we are able to directly compute the invariant
distribution, whereas in the simulation setting, to obtain the same
statistics we would be required to compute a very long
simulations.

If, on the other hand, we simply wish to use the CMA to compute a
trajectory of the slow variables, then the savings will vary, based on
the size of the chosen domain, and the relative differences in
propensity of the fast and slow reactions in the relevant regions. If
our aim is only to produce one relatively short trajectory, then it
is possible that stochastic simulation will be more efficient than
using the CMA. However this is such a trivial task, that any modeller
wishing to do so what not consider invoking any approximations such as
the QSSA or CMA.}

\blue{There are many avenues for future work in this direction, not least
its application to more complex biologically relevant systems. In
particular,  the treatment of systems where the effective behaviour
of the slow variable(s) cannot be well approximated by a
one-dimensional Markov process need to considered, for example systems
which exhibit oscillations. Automated detection of appropriate fast
and slow variables, and statistical tests for the validity of the
approximation for different systems would be hugely beneficial. In the
case of constrained systems which are deficiency zero and weakly
reversible, using the results of \cite{anderson2010product} we can
find the invariant distributions without even constructing the
generator, and this could be a good direction to investigate.}

{\bf Acknowledgements:} The author would like to thank Kostas
Zygalakis for useful conversations regarding this work. This work was funded by First Grant Award EP/L023989/1 from EPSRC.

\bibliography{refs}
\bibliographystyle{plain}
\end{document}